\documentclass[12pt,draftcls,onecolumn]{IEEEtran}
\IEEEoverridecommandlockouts                        

    \usepackage{graphicx}
   \usepackage{epstopdf}
    \usepackage{graphicx}

\usepackage{wrapfig} 
\usepackage{dsfont} 
\usepackage[dvipsnames,usenames]{color}
\usepackage[colorlinks,%
linkcolor=BrickRed,%
filecolor =BrickRed,%
citecolor=RoyalPurple,%
]{hyperref}

\newcommand{\pmat}[1]{\begin{pmatrix} #1 \end{pmatrix}}
\usepackage{amsmath} 
\usepackage{amssymb}  
\usepackage{color}
\usepackage{cite}
\usepackage{verbatim}

\newtheorem{theorem}{Theorem}
\newtheorem{lemma}{Lemma}

\newtheorem{assumption}{Assumption}
\newtheorem{proposition}{Proposition}
\newtheorem{corollary}{Corollary}

\usepackage[loose,normalsize,hang]{subfigure}

\newcommand{\us}[1]{{\color{black}#1}}

 \newcommand{\remove}[1]{}
\newcommand{\EXP}[1]{\mathbb{E}\!\left[#1\right] }

\newcommand{\fhat}[1]{\hat{f} }

\def\Real{\mathbb{R}}

\def\argmin{\mathop{\rm argmin}}

\def\thetahat{{\widehat \theta}}
\def\Fscr{{\mathcal F}}
\def\bhat{{\hat{b}}}

\def\vthetabar{\bar{\vartheta}}
\def\vtheta{\vartheta}
\def\pt{\widetilde{p}}
\def\bbar{\skew3\bar b}

\def\thetabar{\bar\theta}
		\def\qed{\quad \vrule height7.5pt width4.17pt depth0pt}

\usepackage{multirow}

\usepackage{ulem}
\def\use#1{{\color{black}{#1}}}
\def\usd#1{{\color{black}{#1}}}
\def\us#1{{\color{black}{#1}}}
\def\uss#1{{\color{black}{#1}}}
\def\usr#1{{\color{black}{#1}}}

\begin{document}
\title{Distributed computation of equilibria in
			misspecified convex stochastic Nash games}

\author{Hao Jiang \and Uday V. Shanbhag\thanks{\tiny 
Authors may be contacted at \{jh\_double12@hotmail.com,\, udaybag@psu.edu
(corrresponding author),\,
meyn@ece.ufl.edu\} and have  been supported by NSF CAREER CMMI 1246887
	and CMMI 1400217 (Shanbhag). Portions of section \ref{sec:br} are
	based on \cite{jiang11learning}, respectively.
} \and
Sean P. Meyn}
 \maketitle
\vspace{-0.8in}
\begin{abstract}
The distributed computation of Nash equilibria is assuming growing
relevance in engineering where such problems emerge in the
context of distributed control. Accordingly, we present schemes for
computing equilibria of two classes of static
stochastic convex games complicated by a parametric misspecification, a
natural concern in the control of large-scale networked engineered
system. In both schemes,  players {\em learn the equilibrium
strategy} while {\em resolving the misspecification}:
\noindent{(1) \bf Monotone stochastic Nash games:} {We  present a
set of coupled stochastic approximation  schemes distributed across
	agents in which the first scheme updates each agent's strategy via a
	projected (stochastic) gradient step while the second scheme
	updates every agent's belief regarding its misspecified parameter
	using an independently specified learning problem.
	We proceed to show that the produced sequences converge in an
	almost-sure sense  to the true
	equilibrium strategy and the true parameter , respectively.  Surprisingly, convergence in the equilibrium strategy
	achieves the {\em optimal} rate of convergence in a mean-squared sense
	with a quantifiable degradation in the rate constant};
\noindent {(2) \bf Stochastic Nash-Cournot games with unobservable
	aggregate output:} We refine (1) to a Cournot setting where we
	assume that  the
	tuple of strategies is unobservable while payoff
	functions and strategy sets are public knowledge through a common knowledge
			assumption. By utilizing observations of noise-corrupted prices,
{iterative fixed-point} schemes are developed, allowing for {\em
	simultaneously} learning the equilibrium strategies and the
	misspecified parameter in an almost-sure sense.
	\end{abstract}

\section{Introduction}
In  networked engineered systems, a common challenge lies in designing
distributed control architectures that ensure the satisfaction of a
system-wide criterion in environments complicated by nonlinearity,
uncertainty, and dynamics. Such control-theoretic problems take on a
variety of forms and arise in a variety of networked settings, including
networks  of unmanned aerial vehicles (UAVs), traffic networks, wireline
and wireless communication networks, and energy systems, amongst others.
These systems are often characterized by the absence of a designated
central entity that either has system-wide control or has access to
global information. Consequently, control is effected through
distributed decision-making and local interactions that
rely on limited information. Game-theoretic approaches represent an
avenue for designing such protocols. Game theory has seen wide
applicability in the social, economic, and engineered sciences in a
largely {\em descriptive} role. There has been immense
recent interest in a {\em prescriptive role}~\cite{marden13game} that
considers {\em designing a game} whose equilibria represent the solutions to the control problem
of interest~\cite{li10designing,li11distributed}; consequently, the
distributed learning of Nash equilibria assumes immediate relevance in
the management of networked systems.  Learning in Nash games has seen much study in the last several
decades~\cite{Fudenberg,young04strategic,hart05adaptive,shamma05dynamic}. In continuous
strategy regimes, convex static games find
significance in engineered systems such as communication
networks
~\cite{basar07control,alpcan03distributed,pan09games,yin09nash2} and signal
processing~\cite{Fpang09, scutari2011joint}.

An oft-used assumption in game-theoretic models requires that player
payoffs are public knowledge, allowing every player to correctly forecast the
choices of his adversaries. As noted by Kirman~\cite{kirman74learning},
a firm's view of the game may be corrupted or {\em misspecified} in at least two distinct ways
in a Cournot setting where firms decide production levels given a price
function: (i) a firm might have a correct description of the price
function but an incorrect estimate of its parameters; and (ii) it may
have an incorrect structure of the price function and incorrectly
conclude that prediction errors are a consequence of misspecified
parameters.  Kirman~\cite{kirman74learning} considered such a learning
process, and showed that by observing true demand, the suggested
learning process guarantees that the firm strategies converge to the
noncooperative Nash equilibrium. Further inspiration may be drawn from
studies by Bischi~\cite{Bischi07,Bischi08},
Szidarovsky~\cite{Szidarovszky04,szid04stable}, amongst
others~\cite{leonard99nonlinear}, where firms competing in a
deterministic Nash-Cournot game learn a parameter of the demand function
while playing the game repeatedly. Note that an inherent assumption of a
low discount rate is imposed that discounts the future effect of any
player's strategies.  Analogous questions of optimization and estimation have also
been studied by Cooper et al.~\cite{cooper06models} who consider a joint
process of forecasting and optimization in a regime where the underlying
model may be erroneous, demonstrating that the resulting revenues can
systematically reduce over time. 

When designing protocols for practical engineered systems, particularly
in the absence of a centralized controller, the associated parameters of
the utility functions may often be  misspecified. For
instance, in power market models that enlist Nash-Cournot
models~\cite{hobbs01linear,Hobbs07nash}, the precise nature of the price
function is assumed to be given. Similarly, the expected capacity or
availability of renewable generation assets is rarely known a priori.
Similarly, when developing distributed protocols for networked UAVs, the
prescribed utility functions may rely on agent-specific information that
can only be learnt {through} observations. Faced by such challenges,
	our goal lies in the development and analysis of general purpose algorithms that combine
computation of Nash equilibria with a learning phase to correct
misspecification.

\noindent {\bf Motivation:} This research is motivated by the absence of
general-purpose distributed schemes with asymptotic convergence and rate
guarantees for learning equilibria in the
face of imperfect information. 
Such problems emerge from stochastic generalizations
of problems arising in communication
networks~\cite{alpcan03distributed,pavel06noncooperative,pan09games,yin09nash2},
	signal processing~\cite{Fpang09, scutari2011joint}, and
	power markets~\cite{hobbs01linear}.  Accordingly, we present two
	distributed learning schemes in which agents {\em learn their Nash
		strategy} while {\em correcting the misspecification} in their
		payoffs:\\
	\noindent{(1) \bf Stochastic gradient schemes for stochastic Nash
		games:} {In Section~\ref{sec:II}, we present a distributed stochastic approximation
		framework in which every agent makes two projected gradient updates: Every agent
		first updates its belief regarding the equilibrium strategy by
		using the sampled gradient of its payoff function and
		subsequently
		updates its belief regarding the misspecified parameter through
		a similar gradient update.  The resulting
		sequence of equilibrium and parameter estimates are shown to converge to
		their true counterparts in an almost sure sense. Notably, we
		show that the mean-squared error of the equilibrium estimates
		converges to zero at the optimal rate \uss{${\cal O}(1/K)$} despite the presence of
		misspecification where $K$ denotes the number of gradient steps.}\\
	\noindent {(2) \bf \us{Iterative fixed-point} schemes for stochastic
		Nash-Cournot games:} In Section~\ref{sec:br}, we consider a
		Cournot regime where aggregate output is unobservable and one
		parameter of the demand function is misspecified. Under
		common-knowledge, agents develop an estimate of
		aggregate output and the misspecified price function parameter
		by observing noisy prices. These estimates allow developing an
		{iterative fixed-point} scheme that produces iterates that
		are shown to converge to the
		Nash-Cournot equilibrium in an almost-sure sense.  Additionally,
		firms learn the true parameter
		in an almost-sure sense. The result can be extended to nonlinear
		price functions.

\noindent {\bf Remark:} We make two remarks before proceeding. (a)
	First, in (1), the learning problem is constructed independently of the
	computational problem through a set of observations while in (2),
	the learning is affected by the computational step (akin to
			multi-armed bandit problems). (b) Second, we comment on the
{\bf sequential} two-stage framework for resolving
	misspecification:
	$$ \mbox{ Step 1.  Learn $\theta^*$} \qquad \mbox{ Step 2. Compute
		$x^*(\theta^*)$,} $$
	where $\theta^*$ is to be learnt and $x^*(\theta^*)$ is
	the (stochastic) Nash equilibrium, given $\theta^*$. Unfortunately,
{\bf a sequential approach} is complicated by several challenges. First, Step
	1. needs to be completed in a finite number of iterations,
	practically impossible for stochastic learning problems. Second,
	premature termination of Step 1. leads to an erroneous estimate
	$\hat \theta$ \usr{resulting in an incorrect Nash equilibrium
	$\hat x$}. In fact, in stochastic regimes, \usr{such avenues do not
		lead to asymptotic convergence and at best provide approximate
			solutions. We observe from preliminary numerics reveal that sequential schemes may
	perform orders of magnitude worse when compared with iterative
	fixed-point schemes (see Table~\ref{seqsim_fx}).} 

	The rest of the paper is organized as follows. In
	Section~\ref{sec:II}, we define and resolve a misspecified stochastic
	Nash game and present a joint set of stochastic approximation
	schemes that collectively allow for learning equilibria and resolving
	misspecification.  In Section \ref{sec:br}, we
	develop \us{iterative fixed-point} schemes in Cournot settings where  aggregate
	output is unobservable. 
	Empirical studies and conclusions are provided in Sections
	\ref{sec:numerical} and \ref{sec:conc}, respectively.
	Throughout the paper, $\|x\|$ denotes the Euclidean norm of a
vector $x$, i.e., $\|x\|=\sqrt{x^Tx}$ while  $\Pi_K(u)$ denotes the
Euclidean projection of $u$ onto a set $K$, i.e., $\Pi_K(x)\triangleq
\argmin_{y\in K}\|x-y\|$.  A square matrix $H$ is said to be a
$\mathbf{P}$-matrix if every principal minor of $H$ is positive.
Similarly, $H$ is a $\mathbf{P_0}$-matrix if every principal minor of
$H$ is nonnegative.
\section{Gradient-based schemes for convex games}\label{sec:II}

\subsection{Problem description, assumptions and background}  \label{sec:description}
We consider  an $N-$person stochastic
	Nash game in which the
	$i$th player solves Opt$(x_{-i})$:
\begin{align}  \tag{Opt$(x_{-i})$}
\begin{aligned}
	 \min_{x_i \in K_i} & \quad \uss{f_i(x;\theta^*)} \triangleq \EXP{f_i(x;\theta^*,\uss{\xi})} 
\end{aligned}
\end{align}
where $K_i \subseteq \Real^{n_i}$, $\theta^*
\in \Real^m$, $\xi:\Omega \to \Real^d$ \uss{defined on a probability
	space $(\usr{\Omega},\Fscr_x, \mathbb{P}_x)$}, $n = \sum_{i=1}^N n_i$, and $f_i: \Real^n \times \Real^m \times
\Real^d \to \Real$ is a real-valued
function in $x_i$, \uss{$x_{-i}\triangleq(x_j)_{i\neq j=1}^N$}, \usd{$\theta$}
and $\uss{\xi}$. The associated Nash equilibrium is given by a tuple
$x^*=(x_i^*)_{i=1}^N$ where $ x_i^* \in
\textrm{SOL}(\textrm{Opt}(x_{-i}^*))$ for $i=1, \hdots, N,$
$\textrm{SOL}(\textrm{Opt}(x_{-i}^*))$ denotes the solution of
$\textrm{Opt}(x_{-i})$ and under suitable convexity \usd{and differentiability} requirements \usd{(see (A\ref{assump:convex}) below), by invoking~\cite[Th.~7.46]{Shapiro09lecturesSA}}, $x^*$ is a solution to a stochastic variational
inequality problem VI$(K,F(\cdot;\theta^*))$ where
\begin{align}\label{gen-def-F1}
\begin{aligned}
	  K \triangleq \prod_{i=1}^N K_i \mbox{ and }   F(x;\theta) \triangleq \usr{\pmat{
		F_i(x;\theta)}_{i=1}^N \triangleq} \pmat{
		\EXP{\nabla_{x_i} f_i(x;\theta,\uss{\xi})}}_{i=1}^N,
\end{aligned}
\end{align}
respectively. It may be recalled that VI$(K,F)$ requires an $x\in K$ satisfying
\begin{align} \label{VI:gen-main-1}
    (y-x)^T F(x;\theta^*) \geq 0, \quad\, \textrm{for all } y \in K.
\end{align}
Our overall goal lies in computing equilibria when $\theta^*$ is
unavailable or misspecified \usr{but can be learnt by a possibly
	stochastic learning problem.}
\paragraph{Learning scheme}
 \uss{In this section, we consider the estimation of $\theta^*$ through
	 the solution of a suitably defined stochastic convex learning
problem~\cite{hastie01elements}:
\begin{align} \label{problem_stoch_optim_theta}
     \min_{\theta \in \Theta} \quad & g(\theta) \triangleq
	\mathbb{E}[g(\theta;\eta)],
\end{align}
where $\Theta \subseteq \Real^m$ is a closed and convex set,  $\eta:
\usr{\Lambda}\to \Real^p$ is a random variable defined on a probability space
$(\Lambda,\mathcal{F}_{\theta},\mathbb{P}_{\theta})$, and $g:
\Theta \times \usr{\Real^p} \to \Real$ is a real-valued learning metric
function (such as a regression metric constructed from a set of
		observations). Consequently, $\theta^*$ may be learnt through a
stochastic gradient scheme of the form for $k \geq 0$:}
\begin{align}
	\theta_{i}^{k+1}:= \Pi_{\Theta} \left(\theta_i^k - \alpha_{i}^k
			\nabla_{\theta} g(\theta_i^k;\us{\eta_i^k})\right),  \,
	i = 1, \hdots, N.
\end{align}
\us{We emphasize that this learning problem is unrelated to the
computational process and is a built from a set of independently
collected observations.}
\paragraph{Distributed computational scheme} We consider a distributed
stochastic approximation scheme where the $i$th agent employs its
belief regarding $\theta^*$ to take a (stochastic) gradient step for $i = 1, \hdots, N$:
\begin{align}
	x_{i}^{k+1}:= \Pi_{K_i} \left(x_i^k - \gamma_{i}^k \nabla_{x_i}
				f_i(x^k;\usr{\theta_i^k,\xi_i^k})\right), \quad k \geq 0,
\end{align}
\uss{where $\gamma_i^k$ and $\nabla_{x_i}
				f_i(x^k; \usr{\theta_i^k,\xi_i^k})$ denotes the
					steplength and sampled gradient  used by player $i$ at
	step $k$.} \uss{While a fully
rational agent would always take a best response step, in stochastic settings,
the complexity of this step might be significant. In bounded
rational regimes where computational constraints are imposed, an alternative
lies in computing other steps \usr{such as the gradient-response
(cf.~\cite{Papadimitriou94onbounded,Simon96sciences}).}
An \usr{alternative} motivation arises from
distributed control/optimization settings where a ``game'' is designed
whose equilibrium is a desirable solution to a suitably defined control
problem.  Here, a distributed protocol for computing an
equilibrium can be designed and gradient-based approaches can be adopted
(cf.~\cite{marden13game,li10designing,li11distributed}). We propose a game-theoretic
extension of that developed in~\cite{jiang13solution,jiang16solution}.}
We may specify our joint \uss{simulation-based} scheme for learning and computation as follows:
\begin{center}
\fbox{\scriptsize
\parbox[c]{0.98\textwidth}
{{\bf Algorithm I: Gradient response and
	learning}.\\
Let  ${\theta_i^0}\,\in\, \Theta$, $\usr{x^0}\, \in\, K$,
		\uss{$\{\gamma_i^k,\alpha_i^k\}\,>\,0$} ,
	for $i = 1, \hdots, N$, and $k\,=\,0$.\\
{\bf Step 1:}
\uss{\begin{align}\left.
\begin{aligned}
\mathbf{(Comp)}	 \, 	x_{i}^{k+1} & := \Pi_{K_i} \left(x_i^k - \gamma_{i}^k \nabla_{x_i}
				f_i(x^k;{\theta_i^k},\xi_i^k)\right),  \\
	\hspace{-0.2in}\mathbf{(Learn)}	\,  \, \theta_i^{k+1} & := {\Pi_{\Theta }
	\left(\theta_i^k -
			\alpha_i^k \nabla_{\theta } g (\theta_i^k;\eta_i^k)  \right)},
\end{aligned} \right\}  i = 1, \hdots, N;
\end{align}}
{\bf Step 2:} if $k > \bar K$, stop; else $k:=k+1$ and go to Step 1.
}}
\end{center}
\uss{We now present the main assumptions employed in deriving convergence
properties of Algorithm I. (A1) enforces convexity assumptions that allow for
	deriving sufficient equilibrium conditions given by VI$(X,F)$ while the monotonicity
		requirements on $F$ allow for claiming the existence of a unique
		equilibrium. Lipschitzian requirements of $F$ aid in deriving
		subsequent convergence and rate statements. Furthermore, a
		breadth
of learning problems (such as regression, classification
		etc.~\cite{hastie01elements}) are convex. The requirements
imposed by (A\ref{assump:uncommon}) are
standard in developing distributed protocols while (A\ref{assump:filtration})
imposes assumptions on the conditional first and second moments
common in stochastic approximation literature~\cite{Borkar08,Polyak92acceleration,kushner03stochastic}.}

\begin{assumption}[A\ref{assump:convex}] \label{assump:convex}
For $i = 1, \hdots, N$, suppose the function $f_i(x;\theta)$ is convex and
continuously differentiable function in $x_i$ for every $x_{-i}
\in \prod_{j \neq i} K_j$ and every $\theta \in \Theta$. Furthermore,
	suppose $\Theta$ is a closed, convex, and bounded set and for $i = 1,\hdots, N$,
$K_i \subseteq \Real^{n_i}$  is a nonempty, closed, convex and bounded set. Furthermore,  suppose the following hold:(a)  For every $\theta\in\Theta$, $F(x;\theta)$ is both strongly monotone
  and Lipschitz continuous in $x$ with constants $\mu_{x}$ and $L_{x}$; for every $\theta$,
  $(F(x;\theta)-F(y;\theta))^T (x-y)\geq \mu_{x}\|x-y\|^2$ and
  $\|F(x;\theta)-F(y;\theta)\| \leq L_{x}\|x-y\|$;
(b)   For every $x\in K$, $F(x;\theta)$ is Lipschitz continuous in
$\theta$ with constant $L_{\theta}$;
 (c)  The function $g(\theta)$ is strongly convex and continuously differentiable with Lipschitz continuous gradients in $\theta$ with
  convexity constant $\mu_{\theta}$ and Lipschitz constant $C_{\theta}$, respectively\usr{; $(\nabla g(\theta_1)-\nabla g(\theta_2))^T (\theta_1 - \theta_2)\geq \mu_{\theta}\|\theta_1 - \theta_2\|^2$, and
  $\|\nabla g(\theta_1)-\nabla g(\theta_2)\| \leq C_{\theta} \|\theta_1 - \theta_2\|$.}
\end{assumption}

Note that monotone Nash games include {\em stable} Nash
games, a class of games for which a range of
evolutionary dynamics allow for convergence to Nash
equilibria~\cite{hofbauer09stable,fox12population}.
\begin{assumption}[A\ref{assump:uncommon}] \label{assump:uncommon}
For $i=1, \hdots, N$, the $i$th agent knows only his objective $f_i$
and strategy set $K_i$ \usr{and the parameter set $\Theta$}. Furthermore, the vector $x$ is assumed to be
observable.
\end{assumption}

We define a new probability space $(Z,\mathcal{F},\mathbb{P})$, where
$Z\triangleq\Omega\times\Lambda$,
	$\mathcal{F}\triangleq\mathcal{F}_x\times\mathcal{F}_{\theta}$ and
	$\mathbb{P}\triangleq\mathbb{P}_x\times\mathbb{P}_{\theta}$.
	\uss{For $i = 1, \hdots, N$}, suppose $w_i^k \triangleq \nabla_{x_i}
	f_i(x^k;\theta^k_i,\usr{\xi_i^k}) - \nabla_{x_i} f_i(x^k;\theta^k_i)$ and
	$v^k_i \triangleq \nabla_{\theta} g (\theta^{k}_i;\usr{\eta_i^k}) -
	\nabla_{\theta} g (\theta^{k}_i)$.
	$\mathcal{F}_k$ denotes the sigma-field generated by
	$(x^0,\theta^0)$ and errors $(w^l,v^l)$ for $l =
	0,1,\cdots,k-1$, i.e., $\mathcal{F}_0
	=\uss{\sigma}\left\{(x^0,\theta^0)\right\}$ and $ \mathcal{F}_k =
	\uss{\sigma}\left\{(x^0,\theta^0), \left( (w^l,v^l), l = 0,1,\cdots,k-1
			\right)\right\} $ for $k\geq 1.$
	\begin{assumption} [A\ref{assump:filtration}] \label{assump:filtration}
(a) {\em Unbiasedness:} $\mathbb{E}[\usr{w_i^k} \mid \mathcal{F}_k] = 0$ and
  $\mathbb{E}[v_i^k \mid \mathcal{F}_k] = 0$ $a.s.$ for all $k$ and $i$;
  (b) \uss{\em Bounded second moments:} $\mathbb{E}[\|\usr{w_i^{k}}\|^2 \mid \mathcal{F}_k] \leq
	  \nu^2_{x}$ and $\mathbb{E}[\|v_i^{k}\|^2 \mid \mathcal{F}_k]  \leq
		  \nu^2_{\theta}$ $a.s.$ for all $k,i$.
\end{assumption}
 To construct distributed schemes requiring no
coordination in terms of setting parameters, we allow each agent to
independently set steplengths  and as long as a suitable relationship
between these steplengths holds, convergence follows. Specifically, the
$i$th agent employs a diminishing steplength sequence given by
$\gamma_i^k.$  Furthermore, we define $\gamma_{\min}^k\triangleq \min_{1\leq i \leq
N}\{\gamma_i^k\}$ and $\gamma_{\max}^k\triangleq \max_{1\leq i \leq
		N}\{\gamma_i^k\}$ for all $k$. Similarly, we define $\alpha_{\min}^k\triangleq \min_{1\leq i \leq
N}\{\alpha_i^k\}$ and $\alpha_{\max}^k\triangleq \max_{1\leq i \leq
		N}\{\alpha_i^k\}$ for all $k$.  Then, we can make the following assumptions on the steplengths of the algorithm.
\begin{assumption}[Steplength requirements, A\ref{assump:steplength_strongly}] \label{assump:steplength_strongly}
Let $\{\gamma_i^k\}$ and $\{\alpha_i^k\}$ be chosen such that:
(a) $\sum_{k=1}^{\infty}\gamma_{\min}^k = \infty$, $\sum_{k=1}^{\infty}
	(\gamma_{\max}^k)^2 < \infty$, $\sum_{k=1}^{\infty}
	(\alpha_{\max}^k)^2 < \infty$;
(b) $ \lim_{k\to \infty} \frac{\gamma_{\max}^k-\gamma_{\min}^k }{\gamma_{\max}^k} =
	0$; (c)  $\alpha_{\min}^k \geq \gamma_{\max}^k L_{\theta}^2/(\mu_{x}\mu_{\theta})$ for sufficiently large $k$,  $ \lim_{k\to \infty} \frac{(\alpha_{\max}^k)^2}{\gamma_{\max}^k} =
	0.$
\end{assumption}

\uss{Notice that (a)$\sum_{k=1}^{\infty}\gamma_{\min}^k = \infty$ and
	(c)  $\alpha_{\min}^k \geq \gamma_{\max}^k
		L_{\theta}^2/(\mu_{x}\mu_{\theta})$ for sufficiently large $k$
		implies that $\sum_{k=1}^{\infty}\alpha_{\min}^k = \infty$.}
\usr{A natural concern is whether the rule that relates the steplengths can
be implemented in a distributed fashion without coordination. We propose
a rule, first suggested by ~\cite{kannan10online}, in
which  every agent chooses a positive integer and the required
coordination statement holds. We view this as a protocol that may be
employed for developing distributed schemes. The next result ensures
that for such a choice, the required assumptions hold~\cite{kannan10online}.
\begin{lemma}[Choice of steplength sequences] \label{lem:step_choice}
    Let $\{\gamma_i^k\}$ and $\{\alpha_i^k\}$
    be chosen such that    $\gamma_i^k = \frac{1}{(k+N_i)^{\alpha}}$ and $\alpha_i^k = \frac{1}{(k+M_i)^{\beta}}$
	where $N_i$ and $M_i$ are positive integers and $\frac{1}{2} <\beta < \alpha < 1$.  Then,
    \usr{(A\ref{assump:steplength_strongly}) holds.}
\end{lemma}}

We state three results (without proof)
that will be employed in developing our convergence
statements, of which the first two are relatively well-known
super-martingale convergence results (cf.~\cite[Lemma 10,~
Pg.~49--50]{Polyak87})
\begin{lemma}  \label{lem:supermartingale}
Let $\usr{s^k}$ be a sequence of
nonnegative random variables adapted to $\sigma$-algebra $\mathcal{F}_k$ and
such that $\mathbb{E}[\usr{s^{k+1}}\mid \mathcal{F}_k] \leq (1-\usr{u^k})\usr{s^k}+\usr{\beta^k}$ for
  all $ k\geq 0 $ almost surely,
where $0\leq \usr{u^k} \leq 1$, $\usr{\beta^k}\geq 0$, and $\sum_{k=0}^{\infty}\usr{u^k}=\infty$, $\sum_{k=0}^{\infty}\beta^k<\infty$ and $\lim_{k\to\infty}\frac{\beta^k}{u^k}=0$. Then, $\usr{s^k}\to 0$ in an a.s. sense.
\end{lemma}

\begin{lemma}\label{lem:supermartingale2}
Let \usr{$s^k$, $u^k$, $\beta^k$ and $\gamma^k$} be non-negative random
variables adapted to $\sigma$-algebra $\mathcal{F}_k$. \usd{If
\usr{$\mathbb{E}[s^{k+1}\mid \mathcal{F}_k] \leq
  (1+u^k)s^k-\gamma^k+\beta^k$}
   for all  $k\geq 0$ a.s. where
$\sum_{k=0}^{\infty}u^k<\infty$ and $\sum_{k=0}^{\infty}\beta^k<\infty$
both hold a.s., then \usr{$\{s^k\}$} is convergent in an a.s. sense  and
$\sum_{k=0}^{\infty}\gamma_k<\infty$ a.s.}
\end{lemma}
Finally, we present a contraction result reliant on
monotonicity and Lipschitzian requirements.
\begin{lemma}\usd{(cf. \cite[Theorem 12.1.2, Pg. 1109]{Pang03II})} \label{lem:strongly_Lipschitz}
Let $H:K\to\mathbb{R}^n$ be a strongly
		monotone map over $K$ with constant $\mu$, and Lipschitz
		continuous over $K$ with constant $L$. If $q \triangleq
		\sqrt{1-2\mu \gamma + \gamma^2 L^2}$, then for any $\gamma > 0$,
		we have the following:
$\| \Pi_K(x-\gamma H(x))- \Pi_K(y-\gamma H(y))\| \usr{\leq \| x- y - \gamma (H(x) - H(y))\|} \leq q \|x-y\|.$
\end{lemma}

\subsection{Convergence Analysis}

We begin with a contraction statement for the sequence of
iterates produced by Algorithm  I.
\begin{lemma}\label{cexp-bound}
Suppose (A\ref{assump:convex}), (A\ref{assump:uncommon}),
		 (A\ref{assump:filtration})
	and  (A\ref{assump:steplength_strongly}) hold.  Let
	$\{x^k,\theta^k\}$ be computed via Algorithm I. For any $k
	\ge 0$,
$\mathbb{E}\left[\|x^{k+1}-x^*\|^2\mid \mathcal{F}_k\right]  \leq
\zeta_k  \|x^k-x^*\|^2
	 +\beta_k$,
where $\zeta_k  =  1- \gamma_{\max}^k\mu_x + 2 (\gamma_{\max}^k-\gamma_{\min}^k)  L_x +
			2(\gamma_{\max}^k)^2 L_x^2$ and $\beta_k =   ( 2(\gamma_{\max}^k)^2 L_{\theta}^2
    + \gamma_{\max}^k  L_{\theta}^2/\mu_x ) \sum_{i=1}^N \|\theta_i^k-\theta^*\|^2
	    + \usr{N}(\gamma_{\max}^k)^2 \nu_x^2$.
\end{lemma}

{We may now prove our main a.s. convergence result for the sequences
	$\{x^k\}$ and $\{\theta^k\}$.}
\begin{theorem} \label{thm:gradient_strongly_convergence}
Suppose (A\ref{assump:convex}), (A\ref{assump:uncommon}),
		(A\ref{assump:filtration})
	and  (A\ref{assump:steplength_strongly}) hold.  Let
	$\{x^k,\theta^k\}$ be computed via Algorithm I.
Then, $x^{k}\buildrel a.s. \over \longrightarrow x^{*}$ and $\theta_i^{k} \buildrel a.s. \over \longrightarrow   \theta^{*}$ as $k\rightarrow\infty$  for all $i$.
\end{theorem}

Finally, we conclude this section with a non-asymptotic error bound that
demonstrates that the joint scheme displays the optimal rate of
${\cal O}(1/K)$ in mean-squared error.
\begin{theorem} \label{thm:stoch_strongly_optim_error}
\uss{Suppose (A\ref{assump:convex}), (A\ref{assump:uncommon}) and
		(A\ref{assump:filtration}) hold.} Suppose $\gamma_i^k=\lambda_{x,i}/k$
and $\alpha_i^k=\lambda_{\theta,i}/k$.  Let
$\mathbb{E}[\|F_i(x^k;{\theta_i^k}) + w_i^k \|^2]\leq \us{M_x}^2/N$  and
$\mathbb{E}[\|\nabla_{\theta} g(\theta_i^{k})+ v_i^k \|^2]\leq M_{\theta}^2$  for all
$x^k\in K$ and $\theta_i^k\in\Theta$.
Let $\{x^k,\theta^k\}$
be computed via Algorithm I.
Then, \usd{there exist constants $Q_{\theta}$ and $Q_{x,\theta}$} such that the
following hold after $K$ iterations:
  $\mathbb{E}[\|\theta_i^{K} - \theta^{*} \|^2 ] \leq
  Q_{\theta}/K \mbox{ and }
   \mathbb{E}[\|x^{K} - x^{*} \|^2 ] \leq
   Q_{x,\theta}/K.$
\end{theorem}

{\noindent {\bf Remark:} Surprisingly, misspecification does not
lead to a degeneration \usr{in the rate of convergence of the
mean-squared error with respect to that for perfectly specified stochastic Nash games (cf.~\cite{juditsky2011})} but does lead to a worsening of the constant. In
addition, the lack of consistency across steplengths leads to a further
growth in this constant.
}

\section{Iterative fixed-point schemes for misspecified Nash-Cournot games}  \label{sec:br}
\uss{Inspired by the analysis of misspecified Nash-Cournot
	games~\cite{Bischi07,Bischi10,Bischi08,szid04stable,Szidarovszky04},
		we develop an iterative fixed-point scheme. We  introduce the problem in
		Section~\ref{sec:3_description} and describe and analyze the algorithm in
		Sections~\ref{sec_3:desc} and~\ref{sec_3:anal}, respectively. A
		comparison between gradient and iterative fixed-point schemes is
		provided in Section~\ref{sec_3:comp} and we
		conclude with an extension to nonlinear
		prices in~ Section~\ref{sec_3:ext}.}
\vspace{-0.2in}
\subsection{Problem description, assumptions and background} \label{sec:3_description}
We consider a Nash-Cournot game wherein
$ f_i(x) \triangleq c_i(x_i) - p(X;a^*,b^*) x_i,$
where $X \triangleq \sum_{i=1}^N x_i$,
$x_i$ and $c_i(x_i)$ \usr{denote} the scalar output and cost
function associated with  firm $i$ while $K_i$ denotes the strategy set of firm $i$.
Suppose the price function $p(X;a^*,b^*)$ is
defined as
\begin{align}\label{def:price-gen}
  p(X;a^*,b^*)\triangleq\left(a^* -  b^*
			X\right),
\end{align}
Note that $a^*$ represents the ``choke price" at which
demand plummets to zero, while $b^*$ represents the price elasticity
	of demand. Inspired by ~\cite{Bischi08,Bischi10}, we assume that either $a^*$ or $b^*$ is
	unknown  and firm $i$'s belief of this unknown  parameter is denoted
	by $\theta_i$. \usr{We also define $\theta^*$ as the true value of the misspecified parameter of
the price function.} A natural extension is where both parameters are
	unknown and this will require two or more observations at each
	epoch, rather than a single observation of noisy prices. 

\noindent {\bf Case 1 (Learning $a^*$):} We assume that firms know $b^*$ but are unaware
of $a^*$ ($\theta^*=a^*$);  the $i$th firm harbors a belief on $a^*$ denoted
by $\theta_i$ and estimates the aggregate output $X$ by $X_i$, then the $i$th
firm's price estimate and the true noise-corrupted prices are defined as
follows:\begin{align}\label{add-price}
\begin{aligned}
& {p}(X_i;\theta_i,0) \triangleq \theta_i - b^* X_i, && (\mbox{\bf Estimate})\\
& p(X;\theta^*,\xi) \triangleq (\theta^*+\xi) - b^* X.
		&& (\mbox{\bf True price})
\end{aligned}
\end{align}
\noindent {\bf Case 2 (Learning $b^*$):} Distinct from Case 1,  firms
	know $a^*$ and estimate $b^*$ as $\theta_i$ ($\theta^*
			= b^*$) while the
	true price is corrupted by  noise scaled by the aggregate output.
Firm $i$'s price estimate and the true prices are defined as follows:
\begin{align}\label{mult-price}
\begin{aligned}
& {p}(X_i;\theta_i,0) \triangleq a^* - \theta_i X_i, && (\mbox{\bf Estimate})\\
& p(X;\theta^*,\xi) \triangleq a^* - (\theta^* +\xi) X.
		&& (\mbox{\bf True price})
\end{aligned}
\end{align}
The next assumption formalizes these two cases.
\begin{assumption}[A\ref{assump:ab}]  \label{assump:ab}
Either (A\ref{assump:ab}a) or (A\ref{assump:ab}b) holds:\\
    (A\ref{assump:ab}a)  Firms know $b^*$ but not $a^*$ ($\theta^*=a^*$) and the price is
	defined by~\eqref{add-price}. \\
   (A\ref{assump:ab}b) Firms know $a^*$ but not $b^*$ ($\theta^*=b^*$) and the price is defined by
	   \eqref{mult-price}. \\
	   Furthermore, the random variable $\xi$ is defined by
	   $\xi:\Lambda\to \Real,$ $(\Lambda, \Fscr_{\theta}, \mathbb{P}_{\theta})$ is the
	   associated probability space and  $\xi^1,\ldots,\xi^k$ are i.i.d.
	   random variables \uss{with mean zero} for all $k$.
\end{assumption}
Our
assumption on costs is a special case of (A\ref{assump:convex}).

\begin{assumption}[A\ref{assump:convex_cost}] \label{assump:convex_cost}
The cost function $c_i(x_i)$ is a convex and
continuously differentiable function in $x_i$ over $K_i$ with  Lipschitz continuous
gradients with constant
$M_i$. 
Furthermore, $K_1, \hdots, K_N$ are closed, convex, and bounded sets.
\usr{Suppose  the estimator set  $\Theta$ is a compact convex set in $\mathbb{R}_+$
given by $[\delta, \Delta]$ and $0<\delta<\theta^*+\usr{\xi^k}<\Delta$ for all
$k$.}
\end{assumption}

 As forwarded by \cite{aumann76agreeing}, the notion of
``common knowledge'' in game theory extends beyond agents having access to
information. We assume that firms cannot observe
aggregate output and \us{derive an estimate}, relying
on the knowledge of the cost functions and strategy sets of their
competitors, as  assured through a  {\em common
knowledge} assumption. \us{This assumption} is often employed in games (see~\cite{Fudenberg}). Collectively, these two assumptions are captured by (A\ref{assump:common}).
\begin{assumption} [A\ref{assump:common}] \label{assump:common}
The common knowledge assumption holds with regard to $c_i(x_i)$ and $K_i$
for $i=1, \hdots, N.$ Furthermore, aggregate output is unobservable.
\end{assumption}
 Several motivating examples exist in the literature detailing
common knowledge; these include instances provided by
\cite{littlewood53mathematical} (the barbecue problem) and
\cite{schelling60strategy} (the department store problem), amongst
others. While our results are agnostic to
	applications, it is worth emphasizing that such assumptions often
	hold when agents need to make their assets and costs public through
	suitable filings, such as in utility-based regulation (power, gas,
			water, etc.). This is often the case in regulatory settings
	(cf.~\cite[Pg.~78-79]{crew06}). Common knowledge
assumptions immediately hold when a game is
designed~\cite{marden13game,li10designing,li11distributed} and agents can be
endowed with the requisite knowledge. A select number of results will rely on boundedness of strategy
sets, as specified by \usr{(A\ref{assump:convex_cost})}.

\subsection{Description of algorithm} \label{sec_3:desc}
Our goal lies in developing schemes for learning equilibria
and misspecified parameters. Unfortunately, since neither the aggregate
output nor $\theta^*$ are observable, gradient/best-response schemes cannot be directly
implemented. However, under  (A\ref{assump:common}), every firm knows
the cost functions and strategy sets of its competitors, allowing for
computing the best response of all firms, based on an estimate
of $\theta^*$ and the aggregate. By using the discrepancy between
estimated and observed prices, each firm may construct improved
estimates of the misspecified parameter. This model, while aligned, with that suggested by
~\cite{Bischi08,Bischi10} enjoys distinctions at several levels;
specifically, we allow for {\em constrained} problems with {\em
	nonlinear} cost functions with {\em noisy} price observations
	arising from possibly {\em nonlinear} price functions.
	Let $\usd{\bf x}^k_{i}=(x^k_{i1},\cdots,x^k_{iN})$ for
	$i=1,\cdots,N$ and  $X^k_i = \sum_{j=1}^\use{N} x_{ij}^k$  where $x^k_{ij}$ denotes firm $i$'s conjecture of
firm $j$'s output at the $k$th period and $X_i^k$ denote firm $i$'s
estimate of aggregate output.  Note that $X_i^k$ is
maintained as strictly positive by assuming that at least one of the
strategy sets requires strictly positive output while the true aggregate $X^k$ is given
by $X^k \triangleq \sum_{j=1}^N x_{jj}^k$. The proposed algorithm relies on {\bf
	simultaneous}
updates of $x_i^{k+1}$ and $\theta_i^{k+1}$, \usr{which requires first
	defining  $\vtheta_i^{k}$
	$\bar{\vartheta}_i^{k}$, and $\thetahat_i^k$. }
\begin{align}\label{vt-a}
\mbox{under}~\textbf{(A\ref{assump:ab}a)}:
  p(X^{k};\uss{\theta^*},\xi^{k}) & =
				(\theta^*+\xi^{k}) - b^*X^{k}, \quad \vtheta_i^{k}
				\triangleq
				p(X^{\uss{k}};\uss{\theta^*},\xi^{\uss{k}})+b^*X_i^{k},
				 \\ \label{vt-b}
\mbox{under}~\textbf{(A\ref{assump:ab}b)}:
   p(X^{k};\uss{\theta^*},\xi^{k}) & = a^* -
  (\theta^*+\xi^{k})X^{k}, \quad
   \vtheta_i^{k} \triangleq
				(a^*-p(X^{\uss{k}};\uss{\theta^*},\xi^{\uss{k}}))/X_i^{k},
				\\
\label{eq:vtheta}
\usd{\bar{\vartheta}_i^{k+1}}& \usd{=\frac{k\vthetabar_i^{k}+\vtheta_i^{k+1}}{k+1}, } \\
\label{def-hattheta}\widehat{\theta}_i^{k+1}(\theta_i^{k+1}) & \triangleq
\frac{1}{k+1}\theta_i^{k+1}+\frac{k}{k+1} \vthetabar_i^k.
\end{align}
Subsequently, we  show that $\bar{\vartheta}_i^k = \theta^* +
(\sum_{j=1}^k\xi^j)/k$ \usr{(see \eqref{vtheta-savg})}.

\noindent {\bf (a) Update of $x_{i1}^{k+1}, \hdots, x_{iN}^{k+1}$:} Under (A\ref{assump:common}), \usr{given a sequence $\{\epsilon^k\}\downarrow 0$},
firm $i$ \usr{computes a
Nash equilibrium, contingent on its choice
of $\theta_i^{k+1}$, through} a fixed-point of the best-response map:
\begin{align}
   x_{ij}^{k+1} \use{=} \argmin_{\use{x_{ij}^{k+1}} \in K_j}  \ [c_j(x_{ij}^{k+1}) &-
 p(X_{i}^{k+1};\thetahat_i^{k+1}(\theta_i^{k+1})\uss{,0})x_{ij}^{k+1}  +
\uss{\frac{1}{2}}\epsilon^k \|x_{ij}^{k+1}\|^2],
	\quad j = 1, \hdots, N.
\tag{{\bf BR}$_{ij}^x (x_{i,-j}^{k+1},\theta_i^{k+1};
		\usd{\bar{\vartheta}_i^k})$}
\end{align}
\noindent {\bf (b) Update of $\theta_i^{k+1}$:} Firm $i$ defines the
difference between the price
observed at the $k$th step $p(X^{k};\uss{\theta^*},\xi^{k})$ and its
	 estimate ${p}(X^{k+1}_i;\thetahat_i^{k+1},0)$ as
 $\pt_i^{k+1}(\theta_i^{k+1},X_i^{k+1})$:
 \begin{align*}
    \pt_i^{k+1}(\theta_i^{k+1},X_i^{k+1})  := \begin{cases}
				p(X_i^{k+1};\hat{\theta}_i^{k+1}(\theta_i^{k+1}),0) -
				p(X^k;\theta^*,\xi^k), \quad {\scriptsize  \mbox{under }~
				\textbf{(A\ref{assump:ab}a)}}  \\
				p(X^k;\theta^*,\xi^k) -
				p(X_i^{k+1};\hat{\theta}_i^{k+1}(\theta_i^{k+1}),0).
				\quad {\scriptsize  \mbox{under }~
				\textbf{(A\ref{assump:ab}b)}} \end{cases}
\end{align*}
Suppose $t_i^{k+1}(X_i^{k+1})$ denotes a solution to
$$\pt_i^{k+1}(\theta_i^{k+1},X_i^{k+1}) + \epsilon^k \theta_i^{k+1} =0,$$
\use{where the regularization
		term $\epsilon^k \theta_i^{k+1}$ is introduced to ensure uniqueness of
			\eqref{best_response_rand} (See
					Prop.~\ref{lem:best_response_rand}).}
\usr{By invoking the functional form of
	$p(X_i^{k+1};\hat{\theta}_i^{k+1}(\theta_i^{k+1}),0)$,
the following holds:}
\begin{align}\label{eq:t_i}
    t_i^{k+1}(X_i^{k+1})  = \begin{cases}
				[(k+1)(p(X^k;\theta^*,\xi^k) + b^* X_i^{k+1}) - k
				\vthetabar_i^k]/(1+(k+1)\epsilon^k), & \hspace{-0.1in}
{ \scriptsize \mbox{under }~
				\textbf{(A\ref{assump:ab}a)}} \\
				[(k+1)(a^* - p(X^k;\theta^*,\xi^k))- k
				\vthetabar_i^kX_i^{k+1}]/(X_i^{k+1}+(k+1)\epsilon^k).
				&\hspace{-0.1in}
 { \scriptsize \mbox{under }~
				\textbf{(A\ref{assump:ab}b)}} \end{cases}
\end{align}
Suppose $\delta$ and $\Delta$ are lower and upper bounds of $\Theta$, respectively.
We may then update $\theta_i^{k+1}$ as follows:
\begin{align}
\tag{{\bf BR}$_i^{\theta}(X_i^{k+1};\usd{\bar{\vartheta}_i^k})$}
    \theta_i^{k+1} = \begin{cases}
				\delta, \,\;\;\quad & \mbox{if }~
				t_i^{k+1}(X_i^{k+1}) < \delta  \\
				t_i^{k+1}(X_i^{k+1}), \,\;& \mbox{if }~
				\delta\leq t_i^{k+1}(X_i^{k+1}) \leq \Delta  \\
\Delta. \;\quad & \mbox{if }~
				t_i^{k+1}(X_i^{k+1}) > \Delta \end{cases}
\end{align}
\usd{Equivalently, we may state the above as $\theta_i^{k+1} =
	\Pi_{[\delta,\Delta]}(t_i^{k+1}(X_i^{k+1})).$}
%
%
%
\begin{center}
\fbox{\scriptsize
\parbox[c]{0.98\textwidth}
{{\bf Algorithm II: \us{Iterative fixed-point} and
	learning}.\\
Given a sequence \usd{$\{\epsilon^k\}\downarrow 0$ where $\epsilon^k > 0$ for
all $k$},\, and $\gamma_x$,
		$\gamma_{\theta}$; $k := 0; \sum_{j=1}^N x_{jj}^0 = X^0;   p(X^0;\theta^*,\xi^0) := \uss{a^* - b^*X^0};
\epsilon^0 > 0$; \usd{Player $i$ selects $\vthetabar_i^0$} for $i = 1, \hdots, N$. \\
{\bf Step 1.}
For $i = 1, \hdots, N$, if  $ X_i^{k+1}  = \sum_{j=1}^N x_{ij}^{k+1}$,
	then  $\{x_{i1}^{k+1}, \hdots, x_{iN}^{k+1},
	\theta_i^{k+1}\}$ is a solution to the following system:
\begin{align}  \label{best_response_rand}
\begin{aligned}
x_{ij}^{k+1} & \mbox{ solves } \mbox{\bf BR}_{ij}^x
(x_{i,-j}^{k+1},\theta_i^{k+1};\usd{\bar{\vartheta}_i^k}), \quad j = 1, \hdots, N
\\
\theta_i^{k+1} & \mbox{ solves } \mbox{\bf BR}_i^{\theta}(X_i^{k+1};\usd{\bar{\vartheta}_i^k}).
\end{aligned}
\end{align}
{\bf Step 2.} For $i = 1, \hdots, N$, \usd{$\vtheta_i^{k+1}$ is defined as per
\eqref{vt-a} or \eqref{vt-b} and
$\bar{\vtheta}_i^{k+1}$ is updated
as follows:
\begin{align} \label{eq:vtheta-alg}
\bar{\vartheta}_i^{k+1}=\frac{k\vthetabar_i^{k}+\vtheta_i^{k+1}}{k+1}.
\end{align}}
\\{\bf Step 3.} If $k > \bar K$, stop; else $k:=k+1$ and go to Step 1.
}
}
\end{center}
\subsection{Analysis of noise-corrupted \us{iterative fixed-point}
	schemes}\label{sec_3:anal}
In this subsection, we analyze our iterative
	fixed-point scheme and partition the discussion as follows: (i) First,
we provide a brief discussion as to why the update specified by
	\eqref{best_response_rand} can be succinctly
	captured by the solution to a {\bf single}  variational equality problem; (ii) Second, we provide a brief sketch of the results to follows; and (iii) We provide the convergence theory.

\noindent {\bf (i) Equivalence of \eqref{best_response_rand} to a fixed-point problem:}  First,
		any best response of a convex optimization problem is equivalent to a solution of a suitable variational inequality problem~\cite{Pang03I}:
$$ \left[y^*_{\uss{i}} \ \in \ \argmin_{y_\uss{i} \in {\cal Y}_\uss{i}} \ \uss{d}_\uss{i}(y_\uss{i}) \right]
\Leftrightarrow \ \left[ y^*_\uss{i} \
		\mbox{ solves } \mbox{VI}({\cal Y}_\uss{i}, \nabla_{y_\uss{i}} \uss{d}_\uss{i})\right], $$
where $\uss{d}_\uss{i}$ is a convex function in $y_\uss{i}$ over a convex set ${\cal
	Y}_\uss{i}$.  In fact, given a collection of functions $\uss{d}_i(y_i; y_{-i})$ that are convex in $y_i$ over convex sets ${\cal Y}_i$ for all $y_{-i}$ with \uss{$y_{-i}\triangleq(y_j)_{i\neq j}$}, the coupled  best response is equivalent to the solution of a single variational inequality problem\usd{~\cite{Pang03I}}:
\begin{align*}
\left\{\begin{aligned}
\left[y^*_{1} \in  \argmin_{y_1 \in {\cal Y}_1} \ \uss{d}_1(y_1,y_{-1}) \right]
& \Leftrightarrow  \left[ y^*_1 \
		\mbox{ solves } \mbox{VI}({\cal Y}_1, \nabla_{y_1}
				\uss{d}_1(\cdot,y^*_{-1}))\right] \\
				& \ \vdots  \\
\left[y^*_{N} \in \argmin_{y_N \in {\cal Y}_N} \ \uss{d}_N(y_N,y_{-N}) \right]
&  \Leftrightarrow  \left[ y^*_N \
		\mbox{ solves } \mbox{VI}({\cal Y}_N, \nabla_{y_N}
				\uss{d}_N(\cdot,y^*_{-N}))\right]
\end{aligned} \right\} \Leftrightarrow y^* \mbox{ solves } \mbox{VI}({\cal Y}, F),
\end{align*}
where ${\cal Y} \triangleq \prod_{i=1}^N {\cal Y}_i$ and $F(y) =
(\nabla_{y_i} \uss{d}_i(y_i,y_{-i}))_{i=1}^N$. Finally, any solution to a
variational inequality problem is a fixed point of a suitably defined
problem where $\gamma$ is a positive scalar:
$$ \left[y^* \mbox{ solves } \mbox{VI}({\cal Y}, F)\right]
\Leftrightarrow \left[y^* =
\Pi_{\cal Y} (y^* - \gamma F(y^*))\right].$$
By using this avenue, the problem ({\bf
BR}$^x_{ij}(x_{i,-j}^{k+1},\theta_i^{k+1};\usd{\bar{\vartheta}_i^k})$) is the set of coupled
fixed-point problems:
\begin{align} \label{x-fix-pt}
x_{ij}^{k+1}  = \Pi_{K_j} \left(x_{ij}^{k+1} - \gamma \left(\nabla_{x_{ij}}
		f_j(x_i^{k+1};{\thetahat}_i^{k+1}(\theta_i^{k+1})) + \epsilon^k x_{ij}^{k+1}
		\right) \right), \, j = 1, \hdots, N,
\end{align}
where $f_j(x_i^{k+1};{\thetahat}_i^{k+1}(\theta_i^{k+1}))  =
c_j(x_{ij}^{k+1}) -
p(X_{i}^{k+1};\thetahat_i^{k+1}(\theta_i^{k+1}))x_{ij}^{k+1}.$
Similarly, 
({\bf BR}$^{\theta}_{i}(X_{i}^{k+1};\usd{\bar{\vartheta}_i^k})$) can be stated as the following
fixed-point problem:
\begin{align}\label{fix-pt-theta}
\theta_i^{k+1} & = \Pi_{\Theta} \left(\theta_i^{k+1} - \gamma
	\left( \pt_i^{k+1}(\theta_i^{k+1},X_i^{k+1}) + \epsilon^k \theta_i^{k+1} \right)
	\right).
\end{align}
Before proceeding, we shed some light on this equivalence. Suppose the root of $\pt_i^{k+1}(\theta_i^{k+1},X_i^{k+1}) + \epsilon^k
\theta_i^{k+1} = 0$ is denoted by $t_i^{k+1}$. Then from \eqref{fix-pt-theta}, this implies
that  $t_i^{k+1}  = \Pi_{\Theta} \left(t_i^{k+1}\right)$.
Consequently, if
$t_i^{k+1} \in \Theta \triangleq [\delta, \Delta]$, then $\theta_i^{k+1}
= t_i^{k+1}$ while $\theta_i^{k+1} = \delta (\mbox{ or } \Delta)$, if
$t_i^{k+1} < \delta (\mbox{ or } > \Delta)$. But this is equivalent to  ({\bf
		BR}$^{\theta}_{i}(X_{i}^{k+1};\usd{\bar{\vartheta}_i^k})$).

\uss{We define $z_i^{k+1}\triangleq(x_{i1}^{k+1},\hdots,
		x_{iN}^{k+1},\theta_i^{k+1})$. Then, $z_i^{k+1}$ solves the coupled
fixed-point problem }\eqref{x-fix-pt} --
\eqref{fix-pt-theta} if and only if $z_i^{k+1}$ solves VI$({\cal Z}, F^{k+1} \usr{+ \epsilon^k\mathbf{I}})$ where
$$ {\cal Z} \triangleq \prod_{i=1}^N K_i \times \Theta \mbox{ and }
F^{k+1}(\uss{z_i^{k+1}})
	= \pmat{ \nabla_{x_{i1}}
		f_1(x_i^{k+1};{\thetahat}_i^{k+1}(\theta_i^{k+1})) \\
	\vdots \\
	\nabla_{x_{iN}}
		f_N(x_i^{k+1};{\thetahat}_i^{k+1}(\theta_i^{k+1})) \\
\pt_i^{k+1}(\theta_i^{k+1},X_i^{k+1})
	}.$$

\uss{In summary, the coupled best response scheme \eqref{best_response_rand} is equivalent to the coupled fixed-point problem \eqref{x-fix-pt} --
\eqref{fix-pt-theta}, which is also equivalent to the variational inequality problem VI$({\cal Z}, F^{k+1}\usr{+ \epsilon^k\mathbf{I}})$.}\\
\noindent {\bf (ii) Sketch of results:} We first show that the
coupled best response scheme given by \eqref{best_response_rand} always
admits a unique solution (Prop.~\ref{lem:best_response_rand}).
\usr{Furthermore, each player solves a parametrized form of this system in
which the parametrization is shown to be identical
(See~\eqref{vtheta-savg}).} Theorem \ref{thm:br_sto} shows that the
sequence $\usd{\bf x}^k_i,\thetahat_i^k\} \to \{x^*,\theta^*\}$  as $k \to
\infty$ in an a.s. sense. This proof relies on showing that
$\thetahat_i^k \to \theta^*$ as $k \to \infty$ in an a.s. sense. Then,
	if the solution $\usd{\bf x}_i^{k+1}(\thetahat_i^{k+1})$ is a continuous
	function in $\thetahat_i^{k+1}$ (Prop.~\ref{prop:br}), we may
	conclude that $\lim_{k\to \infty} \usd{\bf x}_i^{k+1}(\thetahat_i^{k+1}) =
	\usd{\bf x}_i^{k+1}(\theta^*) = x^*$, where the last equality follows from
	noting that \usr{the correctly specified Nash game has a unique
	solution.}\\
\noindent {\bf (iii) Convergence theory:}

\uss{\begin{proposition} \label{lem:best_response_rand}
Suppose (A\ref{assump:ab}), (A\ref{assump:convex_cost}) and
	(A\ref{assump:common}) hold. \usd{Suppose  $\epsilon^k > 0$ for $k
	\geq 0$.} \usd{G}iven $p(X^k;\theta^*,\xi^k)$ and $\{\vthetabar_i^k\}_{i=1}^N$,
\usd{u}nder either {(A\ref{assump:ab}a)} or
				 {(A\ref{assump:ab}b)},  the
	solution to~\eqref{best_response_rand} is a singleton.
\end{proposition}}

Having shown that the coupled best response scheme has
a unique solution, we proceed to show a Lipschitzian property on the
solution set of \eqref{best_response_rand} with respect to the
parameter $\theta$ (\uss{Prop. \ref{prop:br}}). Before that, we provide some preliminary results.  The strong monotonicity and Lipschitz continuity of the
mapping $F(x)$ can be easily shown under  (A\ref{assump:convex_cost}).
\begin{lemma}\label{mon-Lip}
	Consider the mapping $F(x)$ defined by~\eqref{gen-def-F1} and suppose
		(A\ref{assump:convex_cost}) holds. Then $F(x)$ is a strongly monotone Lipschitz
		continuous mapping.
\end{lemma}
This allows for claiming the existence and uniqueness of a
Nash-Cournot equilibrium when the price function is affine.
\begin{proposition} \label{prop:unique}
	Consider a Nash-Cournot game in which the $i$th player
	solves (Opt$(x_{-i})$) and the price is determined by
	\eqref{def:price-gen}. Furthermore,  suppose (A\ref{assump:convex_cost}) holds. Then, the
	associated Nash-Cournot game admits a unique equilibrium.
\end{proposition}
Now, we state the Lipschitzian property on the
solution set of \eqref{best_response_rand}.  \usr{This proof is inspired by a
	related result (See \cite[Lemma 3]{dafermos88sensitivity}.)}
	\begin{proposition} \label{prop:br}
		Consider a VI$(K,F(\cdot;\theta))$ where $F(x;\theta)$ is strongly
		monotone in $x$ over $K$ for all $\theta \in \Theta$, Lipschitz
		continuous in $x$ for all $\theta \in \Theta$ and Lipschitz
		continuous in $\theta$ for all $x \in K$. Then, the following hold:
	(a) If $x(\theta)$ denotes
		the solution of VI$(K,F(\cdot;\theta))$, then $x(\theta)$ is Lipschitz continuous in
		$\theta$ for all $\theta \in \Theta.$ (b) Given an $\epsilon > 0$, if $x(\theta,\epsilon)$ denotes the solution of
		VI$(K,F(\cdot;\theta)+\epsilon {\bf I})$, then $x(\theta,\epsilon)$
		is Lipschitz continuous in $\theta$ and $\epsilon$.
	\end{proposition}

\uss{If $F_x^{k+1}(\usd{\bf x}_i^{k+1};\thetahat_i^{k+1})
	= \pmat{ \nabla_{x_{ij}}	f_j(\usd{\bf x}_i^{k+1};{\thetahat}_i^{k+1}) }_{j=1}^N,$ then the solution $\usd{\bf x}_i^{k+1}$ to ({\bf
BR}$^x_{ij}(x_{i,-j}^{k+1},\theta_i^{k+1};\usd{\bar{\vartheta}_i^k})$) solves the variational inequality problem
VI$({\prod_{i=1}^N K_i}, F_x^{k+1}(\cdot;\thetahat_i^{k+1})+\epsilon^k\mathbf{I})$\usd{.}
Based on Lemma \ref{mon-Lip} and Prop. \ref{prop:br}, $\usd{\bf x}^{k+1}_i=\usd{\bf x}^{k+1}_i(\thetahat_i^{k+1},\epsilon^k)$ is a continuous function of $(\thetahat_i^{k+1},\epsilon^k)$.
}

	We may now show that the \us{iterative fixed-point} scheme produces a
	sequence of iterates that converge\usd{s} almost surely to the true
	 equilibrium and allow\usd{s} for learning the true parameter.

\begin{theorem}[Global convergence of \us{iterative fixed-point} scheme] \label{thm:br_sto}
	Suppose (A\ref{assump:ab}), \usr{(A\ref{assump:convex_cost}) and (A\ref{assump:common})} hold.
    Let $\{\usd{\bf x}_i^k,\hat{\theta}_i^k\}$ be computed via Algorithm II for
	 $i=1, \hdots, N$. Then $\hat{\theta}_i^k \to \theta^*$  and
	$\usd{\bf x}_{i}^k \to x^*$ almost surely for $i =1, \hdots, N$, where $x^*$ is a solution of the variational inequality
\eqref{VI:gen-main-1}.
\end{theorem}

\noindent {\bf Remark:} We emphasize that this scheme
requires each agent to effectively solve a suitably defined variational inequality
problem, similar to the centralized problem seen
in~\cite{Bischi08,Bischi10}. Such schemes more closely tied to
best-response schemes than the gradient-based approaches presented in
the previous section. \us{In fact, this particular best-response problem requires solving a variational inequality problem with a ${\bf P}$-mapping (See proof of Prop.~\ref{lem:best_response_rand}), a class of problems that has been studied extensively (cf.~\cite[Ch. 10,11]{Pang03I}).}
\vspace{-0.2in}
\usd{\subsection{Comparison between gradient and fixed-point
	schemes}\label{sec_3:comp} To better understand the distinctions between the two schemes, we recap several differences:\\
\noindent {\bf (i) Assumptions:} Gradient-based schemes can accommodate general
convex objectives but require that the gradient map be strongly montone and
Lipschitz continuous. Furthermore, a general misspecified player objective can
be accommodated under suitable assumptions. Iterative fixed-point schemes
focus on stochastic Nash-Cournot games (rather than general stochastic Nash
games) and require convex costs but allow for either linear or nonlinear
inverse demand functions. Notably, players can learn either
$a^*$ or $b^*$ at any given point while a vector of parameters may be learnt in
the context of gradient-based algorithms for stochastic Nash games.\\
\vspace{-0.0251in}
\noindent {\bf (ii) Nature of algorithm:} Gradient-based schemes require
	that each player computes a projected gradient step in the $x$ and
	$\theta$ space while the iterative best-response schemes require
	each player solves a parametrized variational inequality and updates
	a suitably defined average at every step.\\
\vspace{-0.0251in}
\noindent {\bf (iii) Convergence and rate statements:} The gradient-based schemes allow for claiming a.s. convergence, convergence in mean-squared sense, and derive rate statements. In particular, we show that in such regimes the rate of convergence is optimal and no degradation arises from misspecification. On the other hand, in the context of the iterative best-response schemes, we may show a.s. convergence, while rate statements remain a focus of future research.
\vspace{-0.051in}}
\subsection{Extension to nonlinear price functions} \label{sec_3:ext}
We now consider a generalization to nonlinear prices defined as
follows:
\begin{align}\label{price-nonlinear}
	p(X;\theta^*,\xi) \triangleq \begin{cases}
				a^* - b^* X^{\sigma} + \xi, \\
				a^* - (b^*+\xi) X^{\sigma}.
						\end{cases}
\end{align}
This nonlinear price function has been examined
by~\cite{kannan10online} where a discussion of the strict monotonicity of the
associated mapping is presented (Lemma~\ref{monotone-map}(a)). Specifically, the equilibrium of the
Nash-Cournot game are captured by VI$(K,F)$ where $F(x)$ is defined as
\uss{\begin{align}\label{nonlinear-def-map} F(x) \triangleq \pmat{ c_i'(x_i)
	- (a^* - b^*X^{\sigma}) + \sigma
	b^*X^{\sigma -1} x_i}_{i=1}^N.
\end{align}}
In the next result, the mapping $F(x)$ is strongly monotone for all $x \in K$  if $\nabla F(x)$ is a diagonally
dominant matrix for all $x\in K$.
\begin{lemma}\label{monotone-map}
	Consider the mapping $F(x)$ defined in~\eqref{nonlinear-def-map}.
	Suppose  (A\ref{assump:convex_cost}) holds,
	$N < \frac{3\sigma - 1}{\sigma - 1}$ and \usr{$1 < \sigma \leq 3$}. Then the following hold:
	\begin{enumerate}
	\item[(a)] $F(x)$ is a strictly
	monotone mapping over $K$;
	\item[(b)] Suppose $X\geq \eta$ for some $\eta>0$, then $F(x)$ is a
	strongly monotone mapping over $K$.
\end{enumerate}
\end{lemma}

Directly deriving a Lipschitzian statement on $F(x;\theta)$ in terms of
$\theta$ is not easy when the price function has the prescribed
nonlinear form; instead, by noting that $\nabla F(x)$ is bounded when
$x$ is bounded, allows for proving such a statement. Next, we provide a
corollary of Proposition \ref{prop:br} where such a property is derived.
	\begin{corollary}\label{cor:br}
		Consider a VI$(K,F(\cdot;\theta))$ where $F(x;\theta)$ is strongly
		monotone in $x$ over $K$ for all $\theta \in \Theta$, and Lipschitz
		continuous in $\theta$ for all $x \in K$. Also, there is a
		constant $R>0$, such that $\|\nabla F(x;\theta)\|\leq R$ for all
		$x\in K$ and $\theta \in \Theta$. 
Given an $\epsilon > 0$, if $x(\theta,\epsilon)$ denotes the solution of
		VI$(K,F(\cdot;\theta)+\epsilon {\bf I})$, then $x(\theta,\epsilon)$
		is Lipschitz continuous in $\theta$ and $\epsilon$.
	\end{corollary}

\begin{proposition} \label{prop:br_non}
Suppose {(A\ref{assump:ab}a)} holds. Consider the mapping $F(x)$ defined in~\eqref{nonlinear-def-map} and
	suppose \usr{(A\ref{assump:convex_cost}) holds}. Suppose $X\geq \eta$ for some $\eta>0$ and all $x\in K$, where $X=\sum_{i=1}^N{x_i}$. If
	$N < \frac{3\sigma - 1}{\sigma - 1}$ and \usr{$1 < \sigma \leq 3$}, then the following hold:
		\begin{enumerate}
			\item[(a)] If $x(\theta)$ denotes
		the solution of VI$(K,F(.;\theta))$, then $x(\theta)$ is Lipschitz continuous in
		$\theta$ for all $\theta \in \Theta.$
		\item[(b)] Given an $\epsilon > 0$, if $x(\theta,\epsilon)$ denotes the solution of
		VI$(K,F(.;\theta)+\epsilon {\bf I})$, then $x(\theta,\epsilon)$
		is Lipschitz continuous in $\theta$ and $\epsilon$.
		\end{enumerate}
	\end{proposition}

We may now show that the \us{fixed-point} problem yields a unique solution.
\begin{proposition} \label{prop:br_non_unique}
Suppose (A\ref{assump:convex_cost}) and (A\ref{assump:common})  hold. Let the price be given by
\eqref{price-nonlinear}. If $N < \frac{ 3\sigma - 1}{\sigma - 1}$ and
\usr{$1 < \sigma \leq 3$}, then given $p^k(\xi^k)$ and $\{\thetabar_i^k\}_{i=1}^N$,
	the solution to~\eqref{best_response_rand} is a singleton.
\end{proposition}

By leveraging Propositions \ref{prop:br_non} and \ref{prop:br_non_unique}, the convergence of the
\us{iterative fixed-point} scheme can be claimed under the caveat that
the aggregate output is always bounded away from zero, as stated by the
next result, whose proof is similar to
Theorem~\ref{thm:br_sto} and is omitted.
\begin{corollary}
Suppose \usr{(A\ref{assump:convex_cost}) and (A\ref{assump:common})} hold. Suppose $X\geq \eta$ for some $\eta>0$ and all $x\in K$, where $X=\sum_{i=1}^N{x_i}$.
    Let $\{\usd{\bf x}_i^k,\hat{\theta}_i^k\}$ be computed via Algorithm~II for
	 $i=1, \hdots, N$.
Suppose a unique solution to the \us{fixed-point} problem \eqref{best_response_rand} can be obtained, given $p^k(\xi^k)$ and $\{\thetabar_i^k\}_{i=1}^N$ for each $k\geq 0$.
Then, $\hat{\theta}_i^k \to \theta^*$ almost surely for $i=1, \hdots, N$ and
	$\usd{\bf x}_{i}^k \to x^*$ almost surely for $i =1, \hdots, N$, where $x^*$ is a solution of the variational inequality \eqref{VI:gen-main-1}.
\end{corollary}

We conclude this section with an observation. If one
	used a more widely used estimation technique such as a least-squares
	estimation then it remains unclear if almost-sure convergence
	statements can always be claimed since least-squares
	estimators generally converge in a weaker sense while stronger
	statements may be available for linear regression (see \cite{Anderson79LeastSquares}). In
	effect, a scheme that combines a least-squares estimation technique
	with a strategy update, while convergent, {\em may} not possess desirable
	almost-sure convergence properties.   While, we examine nonlinear
	Nash-Cournot games in this section, we also show that
	such claims hold for more general aggregative Nash games.
	However, it should be emphasized that extending this avenue to Nash
	games where the associated variational map is non-monotone may
	lead to challenges. In particular, what are perfectly reasonable
	schemes for a subclass of Nash games may not be supported by similar
	asymptotic guarantees when the structural properties of the problem
	do not satisfy some key requirements.

\section{Networked Nash-Cournot games}  \label{sec:numerical}
In this section, we apply the developed algorithms
on a class of networked Nash-Cournot games described in Section~\ref{sec:61}.
In Section~\ref{sec:63}, we apply the distributed gradient-based schemes for
purposes of learning equilibria and the misspecified parameters when
aggregate output is observable,
while in
Section~\ref{sec:62}, we apply the proposed \us{iterative fixed-point} schemes  when
aggregate output is unobservable.
Note that the simulations were carried
out on Matlab R2009a on a laptop with Intel Core 2 Duo CPU (2.40GHz) and
2GB memory.  The complementarity solver \texttt{PATH}, developed by
~\cite{Ferris98complementarity}, was utilized for solving the
variational inequality problems that arose in implementing the
algorithms.
\subsubsection{Problem description}\label{sec:61}
We consider a  setting where there are $N$
firms competing over a $W$-node network.  Firm $f$ may produce and sell
its good at node $i$ (denoted by $g_{fi}$ and $s_{fi}$, respectively), where $f=1,\ldots,N$ and
 $i=1,\ldots,W$.  We assume that for a given firm $f$, the cost of
 generating $g_{fi}$ units of power at node $i$ is linear and is given
 by $c_{fi}g_{fi}$. Furthermore, the generation level associated with
 firm $f$ is bounded by its  production capacity,  which is denoted by $\textrm{cap}_{fi}$.
The aggregate sales of all firms at node $i$ is denoted by $S_{i}$, and
the nodal price of power at node $i$, assumed to be a linear function
of $S_{i}$, is defined as
  $  p_i(S_i) \triangleq a^*_i - b^*_i S_i$,
where $a_i^*$ and $b_i^*$ are node-specific positive price function
parameters.  A given firm can produce at any node and then sell at different
nodes, provided that the aggregate production at all nodes matches the
aggregate sales at all nodes for each firm. For simplicity, we assume
that there is no transportation cost between any two nodes, and
that there is no limit of sales at any node.  Then, the resulting
problem faced by firm $f$ can be stated as
\begin{align}
&(\textrm{Firm}(x_{-f}))   \label{prob:nash_network}    \max_{s_{fi} \geq 0,\, \textrm{cap}_{fi} \geq g_{fi} \geq 0} \left\{
	\sum_{i=1}^W \left( p_i(S_i)s_{fi} - c_{fi}g_{fi} \right): \sum_{i=1}^W (s_{fi}-g_{fi}) = 0
		 \right\}.
\end{align}
The resulting Nash-Cournot equilibrium is given by $\{x_f^*\}_{f=1}^N$
where $x_f^*$ is a solution to (Firm$(x_{-f}^*)$) for $f= 1, \hdots, N$.
Prices are assumed to be corrupted by noise, in one of two ways:
\begin{align}  \label{network_price_a}
    p_i(S_i;\xi_i) & = (a^*_i + \xi_i) - b^*_i S_i, \\
  \label{network_price_b}
    p_i(S_i;\xi_i) & = a^*_i - (b^*_i+ \xi_i) S_i.
\end{align}
Note that firm $f$ either has to learn $\theta^*\triangleq
(a_i^*)_{i=1}^W$ when prices are given by \eqref{network_price_a} or
learn $\theta^*\triangleq (b_i^*)_{i=1}^W$ when prices are given by
\eqref{network_price_b}.  In the remainder of this section,  let \break
$a^*\triangleq(a_1^*,\ldots,a_W^*)^T$,
   $b^*\triangleq(b_1^*,\ldots,b_W^*)^T$,
   $\theta^*\triangleq(\theta_1^*,\ldots,\theta_W^*)^T$,
   $\xi^*\triangleq(\xi_1^*,\ldots,\xi_W^*)^T$, $x_i\triangleq(s_{1i}, s_{2i}, \ldots,
		   s_{Ni},  g_{1i}, g_{2i}, \ldots, g_{Ni})^T$, and
   $x\triangleq(x_1^T,\ldots,x_W^T)^T$.  Note that this problem is employed as a motivating example
since Cournot-based models have been used extensively in their analysis
(cf.~\cite{hobbs01linear,Hobbs07nash}). Naturally, a range of
rationality assumptions can be imposed on firms in power markets, but
given the sheer size of the problem and the repeated nature of
competition (in most power markets, firms compete as many as 5--6 times
		every hour in the setting of prices) with relatively minor
changes occuring in demand/availability over a short period.
\subsubsection{Learning with observation of the aggregate
	output}\label{sec:63}
In this subsection, we assume that every firm knows the aggregate output at each node, and employ the learning schemes proposed in Section \ref{sec:description}.  Suppose, the nodal price function is given by (\ref{network_price_a}) and suppose Algorithm I (the
		gradient-based distributed learning scheme), proposed in Section
\ref{sec:description}, is employed for  learning parameters and computing equilibria.
Suppose firms have generated a price at each node. We use $p_{i} = a_{i}^* +\xi_i - b_{i}^*S_{i}$ to denote the
price.
Each firm will solve the following (regularized) problem to estimate $a_{i}^*$ and $b_{i}^*$:
\begin{align*}
 \begin{aligned}
    \min_{ \usr{\{a_{i} ,b_{i}\} \triangleq \theta} \in \Theta}\quad \usr{g(\theta)},  \mbox{ where }
 \usr{g(\theta)} \triangleq   \mathbb{E}\left[  ( a_{i}  - b_{i}S_{i} - p_{i})^2 +
     \lambda  a_{i}^2 +  \lambda
	 b_{i}^2\right].
\end{aligned}
\end{align*}
Suppose $S_{i}$ is as per a uniform distribution
and is specified  by $S_{i}\sim U[0, a_{i}^0/b_{i}^0]$, where $a_{i}^0$ and $b_{i}^0$ are initial estimates of $a^*_i$ and $b^*_i$. Suppose, the noise $\xi_i$ is distributed as per a uniform distribution
and is specified  by $\xi_i\sim U[-{a}^*/2,{a}^*/2]$. Suppose
the steplength sequence $\{\gamma_i^k\}$ and $\{\alpha_i^k\}$ are chosen according to
Lemma \ref{lem:step_choice}: $\gamma_i^k = \frac{1}{(k+N_i)^{\alpha}}$ and $\alpha_i^k = \frac{1}{(k+M_i)^{\beta}}$,
	where $\alpha=0.8$ and $\beta=0.6$ and $N_i$ and $M_i$ are randomly chosen from an interval
	$[1, 200].$ The algorithm was terminated at $k=10000$.  Table
	\ref{table:grad_sto} shows the scaled errors of the learning scheme.

\begin{table}[htbp]
\centering \hspace{-.1in}
\begin{minipage}{0.38\textwidth}
\centering
\caption{Distributed gradient scheme}   \label{table:grad_sto}
\tiny
\vspace{-0.1in}
\begin{tabular}{|c|c|c|c|c|c|c|c|c|c|c|c|}
  \hline
   \multirow{2}{*}{N} & \multirow{2}{*}{W}& \multicolumn{3}{c|}{Learning $a^*$ and $b^*$}  \\
   \cline{3-5} & &  $\frac{\|x^{k}- x^*\|}{1+\|x^*\|}$ & $\frac{\|\hat{a}^k- a^*\|}{1+\|a^*\|}$  & $\frac{\|\hat{b}^k- b^*\|}{1+\|b^*\|}$   \\
 \hline
  5 &  	  1 &	 7.2$\times 10^{-7}$  &	 2.9$\times 10^{-2}$  & 4.7$\times 10^{-2}$  \\
  \hline
  5  &	  2 & 	  3.3$\times 10^{-4}$  &	 3.3$\times 10^{-2}$  & 5.3$\times 10^{-2}$ \\
  \hline
   5 & 	  3 &    7.4$\times 10^{-5}$  &	 3.3$\times 10^{-2}$  & 5.3$\times 10^{-2}$    \\
  \hline
   5 & 	  4 &    1.2$\times 10^{-2}$  &	 4.2$\times 10^{-2}$  & 6.8$\times 10^{-2}$  \\
  \hline
   5 & 	  5 &   1.4$\times 10^{-2}$  &	 3.2$\times 10^{-2}$  & 8.5$\times 10^{-2}$  \\
    \hline
   10 & 	  2 &   1.3$\times 10^{-4}$  &	 3.4$\times 10^{-2}$  & 3.7$\times 10^{-2}$   \\
   \hline
   10 & 	  4 &   1.1$\times 10^{-2}$  &	 2.6$\times 10^{-2}$  & 8.4$\times 10^{-2}$   \\
   \hline
   10 & 	  6 &   2.4$\times 10^{-2}$  &	 3.6$\times 10^{-2}$  & 8.6$\times 10^{-2}$   \\
   \hline
   10 & 	  8 &   2.8$\times 10^{-2}$  &	 3.0$\times 10^{-2}$  & 6.4$\times 10^{-2}$   \\
   \hline
   10 & 	  10 &   3.1$\times 10^{-2}$  &	 4.1$\times 10^{-2}$  & 5.4$\times 10^{-2}$   \\
  \hline
\end{tabular}
\end{minipage}\,
\begin{minipage}{0.61\textwidth}
\centering
\caption{Iterative fixed-point scheme}   \label{table:br_sto}
\tiny
\vspace{-0.1in}
\begin{tabular}{|c|c|c|c|c|c|}
  \hline
  \multirow{2}{*}{N} & \multirow{2}{*}{W}& \multicolumn{2}{c|}{Learning $a^*$} & \multicolumn{2}{c|}{Learning $b^*$} \\
   \cline{3-6} & & $\max_f\frac{\|x_f^{k}- x^*\|}{1+\|x^*\|}$ & $\max_f\frac{\|\hat{a}_f^k- a^*\|}{1+\|a^*\|}$ & $\max_f\frac{\|x_f^{k}- x^*\|}{1+\|x^*\|}$ & $\max_f\frac{\|\hat{b}_f^k- b^*\|}{1+\|b^*\|}$  \\
  \hline
   5 &  	  1 &	  6.0$\times 10^{-3}$ &  5.4$\times 10^{-3}$  &  2.3$\times 10^{-3}$  &	  1.5$\times 10^{-3}$\\
  \hline
  5  &	  2 & 	  1.9$\times 10^{-3}$ &	  1.6$\times 10^{-3}$  & 	   9.1$\times 10^{-4}$  &	  7.7$\times 10^{-4}$ \\
  \hline
   5 & 	  3 &     1.4$\times 10^{-3}$ &  2.7$\times 10^{-3}$ &     7.8$\times 10^{-4}$  &	  1.4$\times 10^{-3}$ \\
  \hline
  5 & 	  4 &     7.8$\times 10^{-3}$ &  2.8$\times 10^{-3}$ &   2.0$\times 10^{-3}$  &	  1.0$\times 10^{-3}$ \\
  \hline
  5 & 	  5 &     1.0$\times 10^{-3}$ &	  2.5$\times 10^{-3}$ & 1.2$\times 10^{-2}$  	&  2.2$\times 10^{-3}$ \\
  \hline
  10 & 	  2 &     2.0$\times 10^{-3}$ &  1.9$\times 10^{-3}$ &	  1.2$\times 10^{-3}$  	&  1.2$\times 10^{-3}$\\
  \hline
  10 & 	  4 &    1.1$\times 10^{-2}$ &	   4.2$\times 10^{-3}$ &    1.5$\times 10^{-2}$  & 	  9.4$\times 10^{-4}$ \\
  \hline
  10 & 	  6 &    1.8$\times 10^{-3}$ &	  0.8$\times 10^{-3}$ &    3.0$\times 10^{-4}$  &	  1.5$\times 10^{-3}$ \\
  \hline
  10 & 	  8 &     2.0$\times 10^{-3}$&	   2.7$\times 10^{-3}$ &    1.3$\times 10^{-3}$    &     8.5$\times 10^{-4}$  \\
  \hline
  10 & 	  10 &    1.1$\times 10^{-3}$&	  3.5$\times 10^{-3}$ &   3.8$\times 10^{-4}$    &      7.0$\times 10^{-4}$\\
  \hline
\end{tabular}
\end{minipage}
\end{table}


\subsubsection{Learning without observing the aggregate
	output}\label{sec:62}
In this subsection, we examine how the schemes perform when firms are
ignorant of aggregate output at each node while a common knowledge
assumption is assumed to hold.

%

\begin{figure}[htb]
\subfigure[Learning $x$]{  \begin{minipage}[b]{0.5\textwidth}
    \centering
    \includegraphics[width=0.8\textwidth]{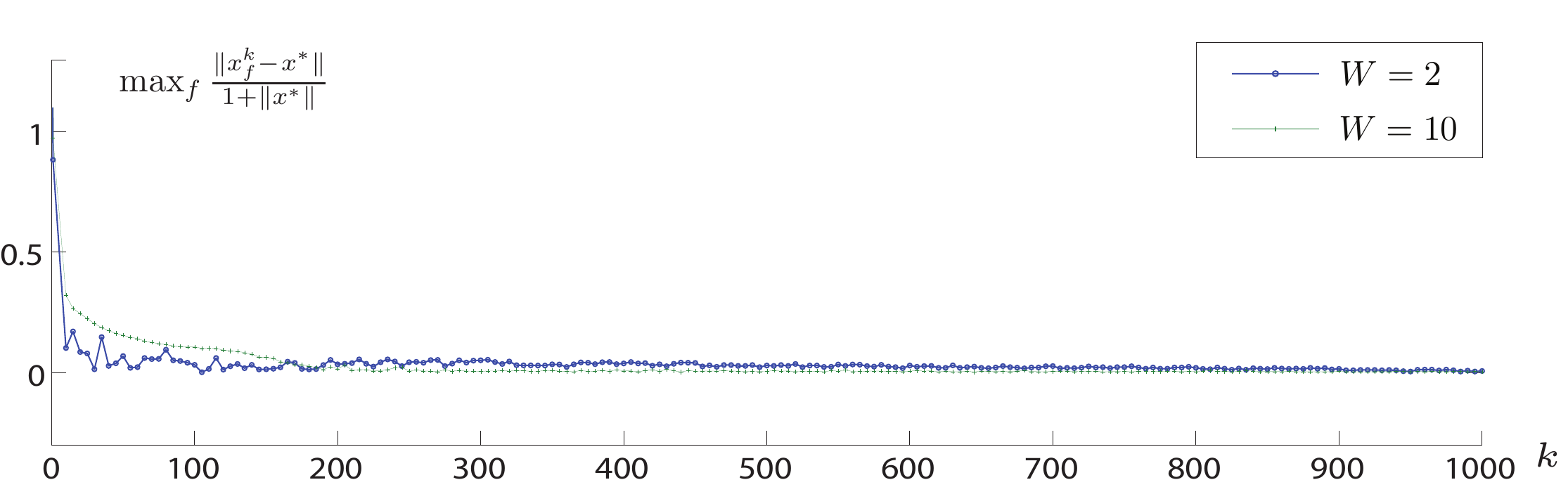}
   \label{figure:br_sto_ax}
  \end{minipage}}
\subfigure[Learning $a$]{  \begin{minipage}[b]{0.5\textwidth}
    \centering
    \includegraphics[width=0.8\textwidth]{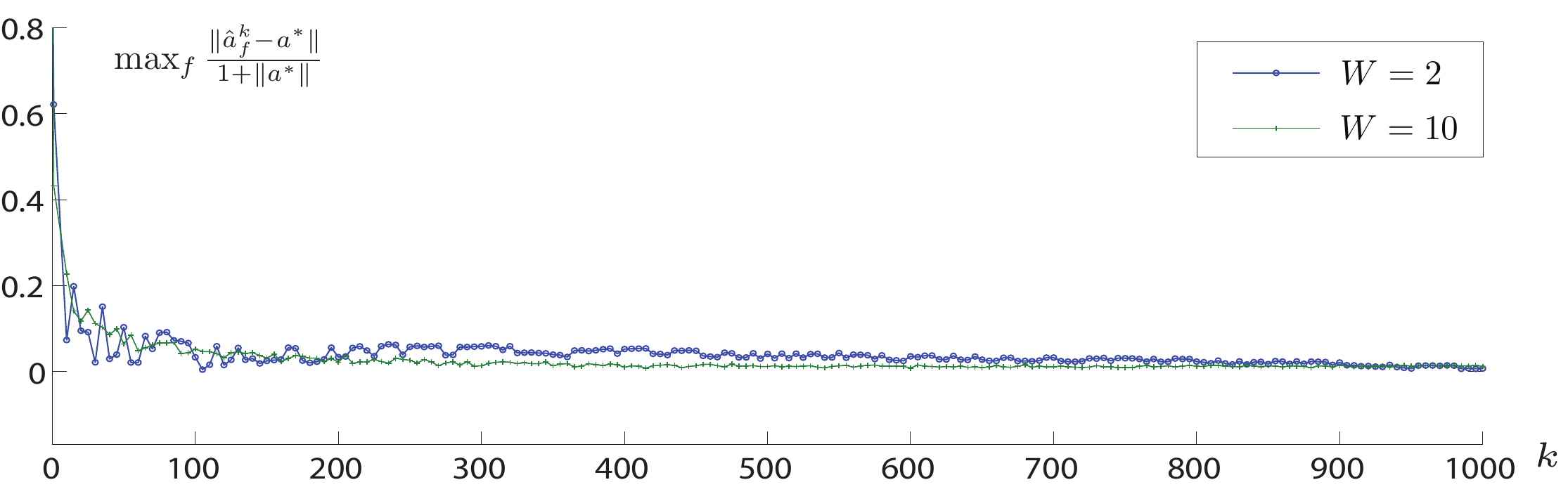}
     \label{figure:br_sto_aa}
  \end{minipage}}
  \caption{Computing $x^*$ and learning $a^*$ ($\xi\sim U[-{\theta}^*/2,{\theta}^*/2]$,
			$N=10$)}
\end{figure}

\begin{figure}[htb]
\subfigure[Learning $x$]{  \begin{minipage}[b]{0.5\textwidth}
    \centering
    \includegraphics[width=0.9\textwidth]{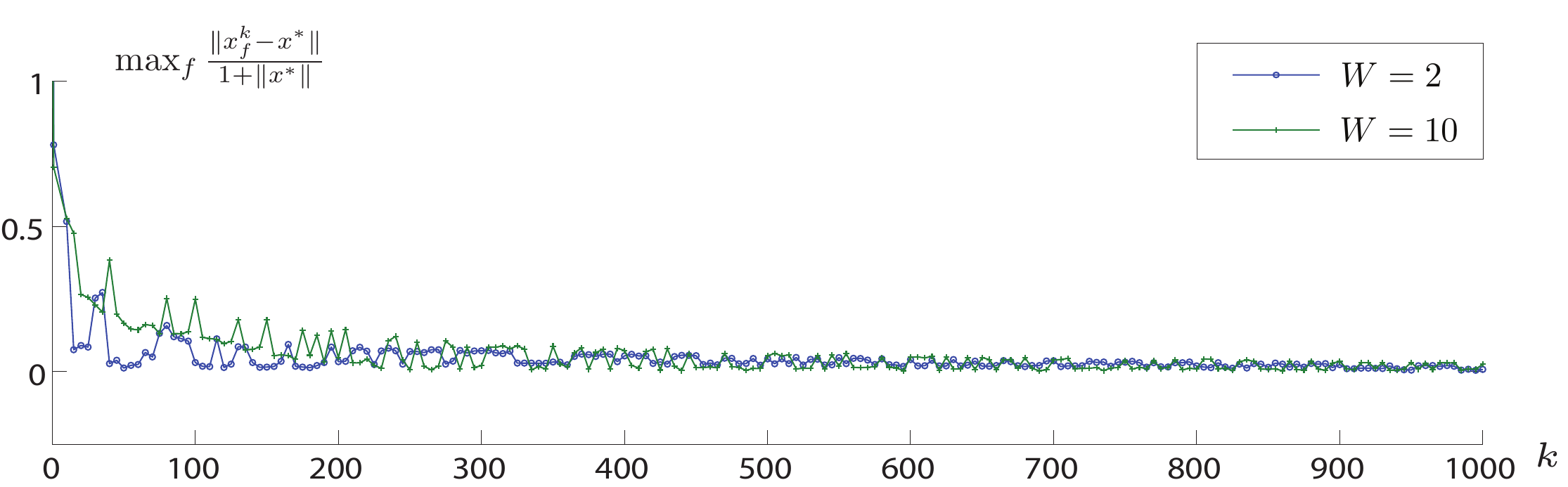}
   \label{figure:br_sto_bx}
  \end{minipage}}
\subfigure[Learning $b$]{  \begin{minipage}[b]{0.5\textwidth}
    \centering
    \includegraphics[width=0.9\textwidth]{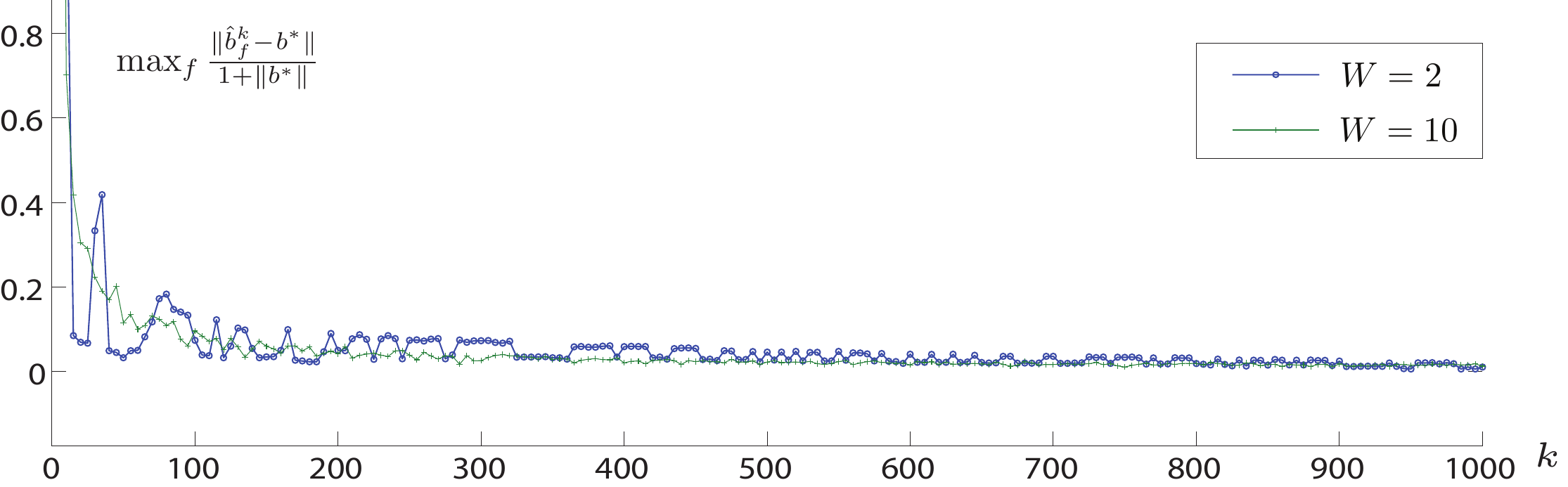}
     \label{figure:br_sto_bb}
  \end{minipage}}
\caption{Computing $x^*$ and learning $b^*$ ( $\xi\sim U[-{\theta}^*/2,{\theta}^*/2]$, $N=10$)}
\end{figure}

Suppose, the nodal price function is given by (\ref{network_price_a}) or
(\ref{network_price_b}) and suppose Algorithm II (the
		\us{iterative fixed-point} scheme), proposed in Section
\ref{sec:3_description}, is employed for  learning parameters and computing equilibria.
Suppose, the noise $\xi$ is distributed as per a uniform distribution
and is specified  by $\xi\sim U[-{\theta}^*/2,{\theta}^*/2]$. Each run
comprised of 10000 steps learning $a^*$ and 50000 steps for learning
$b^*$.  Table \ref{table:br_sto} shows the scaled errors of the learning
scheme while Figures \ref{figure:br_sto_ax} and \ref{figure:br_sto_aa}
illustrate the scaled errors of the learning scheme when the number of
steps, denoted by $k$, increases for learning $x^*$ and $a^*$,
respectively. Analogous figures for learning $x^*$ and $b^*$ are
provided (see Figures \ref{figure:br_sto_bx} and
		\ref{figure:br_sto_bb}).  


\vspace{-.1in}
\begin{table}[htbp]
\caption{Learning $x^*$ and $b^*$ in a stochastic regime when $N=5$ and $W=1$, stopping at step $k=10000$}
\tiny
\vspace{-0.2in}
\centering
\subtable[$\xi\sim U\textrm{[}-b^*/2,b^*/2\textrm{]}$]{
\begin{minipage}{0.49\textwidth}\label{tab1}
\centering\vspace{-0.1in}
\begin{tabular}{|c|c|c|c|c|c|c|c|c|c|c|c|}
  \hline
  & \multicolumn{2}{|c|}{Sequential} & \multicolumn{2}{|c|}{Simultaneous}\\
  \hline
  Bound& $\frac{\|x^{k}- x^*\|}{1+\|x^*\|}$   & $\frac{\|\hat{b}^k- b^*\|}{1+\|b^*\|}$ & $\frac{\max\|x_f^{k}- x^*\|}{1+\|x^*\|}$ & $\frac{\max\|\hat{b}_f^k- b^*\|}{1+\|b^*\|}$ \\  
 \hline
 32.3664 &	 2.1$\times 10^{-1}$  	&  1.2$\times 10^{-1}$ &	 4.9$\times 10^{-3}$  	&  3.3$\times 10^{-3}$ \\
  \hline
   64.7329 &	 1.2$\times 10^{-1}$ & 1.0$\times 10^{-1}$  &	 5.0$\times 10^{-3}$ & 3.3$\times 10^{-3}$  \\
  \hline
   97.0993  &	5.5$\times 10^{-1}$ & 8.8$\times 10^{-1}$   &	 5.0$\times 10^{-3}$ & 3.3$\times 10^{-3}$  \\
  \hline
   129.4658 &	 7.4$\times 10^{-1}$ & 1.1  &	 5.1$\times 10^{-3}$ & 3.4$\times 10^{-3}$ \\
  \hline
   161.8322 &	 1.2 & 7.9$\times 10^{-1}$  &	 5.1$\times 10^{-3}$ & 3.4$\times 10^{-3}$ \\
  \hline
\end{tabular}
\end{minipage}}\hfill
\subtable[$\xi\sim U\textrm{[}-R,R\textrm{]}$]{
\begin{minipage}{0.49\textwidth} \label{tab2}
\centering\vspace{-0.1in}
\begin{tabular}{|c|c|c|c|c|c|c|c|c|c|c|c|}
   \hline
  & \multicolumn{2}{|c|}{Sequential} & \multicolumn{2}{|c|}{Simultaneous}\\
  \hline
  $R$ & $\frac{\|x^{k}- x^*\|}{1+\|x^*\|}$   & $\frac{\|\hat{b}^k- b^*\|}{1+\|b^*\|}$ & $\frac{\max\|x_f^{k}- x^*\|}{1+\|x^*\|}$ & $\frac{\max\|\hat{b}_f^k- b^*\|}{1+\|b^*\|}$ \\   
 \hline
 $b^*/5$ &	7.5$\times 10^{-2}$ & 4.8$\times 10^{-2}$ &	 1.9$\times 10^{-3}$  	&  1.2$\times 10^{-3}$ \\
  \hline
   $b^*/4$ &	9.6$\times 10^{-2}$ & 6.0$\times 10^{-2}$  &	 2.4$\times 10^{-3}$ & 1.6$\times 10^{-3}$  \\
  \hline
   $b^*/3$ &	 1.3$\times 10^{-1}$ & 8.0$\times 10^{-2}$   &	 3.2$\times 10^{-3}$ & 2.2$\times 10^{-3}$  \\
  \hline
   $b^*/2$ &	 2.1$\times 10^{-1}$  	&  1.2$\times 10^{-1}$ &	 4.9$\times 10^{-3}$  	&  3.3$\times 10^{-3}$ \\
  \hline
   $b^*/1$ &	 5.3$\times 10^{-1}$  	&  2.3$\times 10^{-1}$ &	 9.9$\times 10^{-3}$  	&  6.7$\times 10^{-3}$ \\
  \hline
\end{tabular}
\end{minipage}}
\label{seqsim_fx}
\end{table}

\vspace{.3in}

In Table
\ref{tab1}, we raise the upper bounds of the
strategy sets of all agents and compare a sequential scheme with our
iterative fixed-point scheme. \us{In the sequential counterpart, we
	employ $10,000$ steps of stochastic approximation-based learning followed by $10,000$ steps of
		computation.}  It is seen that the error
from the sequential scheme
increases proportionally to the bound, while the error associated with
our simultaneous scheme does not change significantly.  Table \ref{tab2} shows that
when increasing the variance of the noise makes the difference in errors between the
sequential and simultaneous schemes more pronounced.
Consequently, for the same effort, it can be seen that the simultaneous
scheme performs far better to the sequential scheme, particularly when
the variance of the noise grows.

\vspace{-.2in}

\section{Concluding remarks}\label{sec:conc}
Nash games, a broadly applicable paradigm for modeling
strategic interactions in noncooperative settings, have emerged as
immensely useful in the context of distributed control problems.  Yet,
the development of distributed protocols for learning equilibria may be
complicated by several challenges: (i) Agents may have an incomplete
specification of payoffs; (ii) Agents may be unavailable to observe the
actions of their counterparts; and finally, (iii) Observations  may be
corrupted by noise.  Accordingly, this paper is motivated by developing
schemes for learning
equilibria and resolving misspecification (such as in the price
		functions). We consider two specific settings as
part of our investigation and apply these techniques on a class of
networked Nash-Cournot games.  First, we consider convex static
stochastic Nash games characterized by a suitable monotonicity property
in which agent payoffs are parameterized by a  misspecified vector. We
consider a framework that combines (stochastic) gradient steps with a
stochastic approximation step that attempts to learn the parameter.   In such settings, we provide asymptotic
statements that show that
agents	may learn equilibria and the true parameters in an almost sure
sense. In addition, we provide non-asymptotic error bounds that
demonstrate that the rate of convergence is not impaired by the presence
of learning.  Second, we refine our statements to a Cournot regime where we assume {\em common
	knowledge} holds but aggregate output is
	unobservable. In such a setting, we construct a learning scheme in
	which firms maintain a belief of the aggregate output and the
	misspecified price function parameter. After each step, these
	beliefs are updated by employ \us{fixed-point} steps and by leveraging
	the disparity between estimated and (noisy) observed prices. We proceed to
	show that in the limit, every firm learns the true Nash-Cournot equilibrium
	strategy in an almost-sure sense. Additionally, every firm learns
	the correct value of the misspecified parameter in an almost-sure
	sense.  {Yet much remains
		to be studied, including weakening monotonicity
		requirements on the map and boundedness requirements on the
		strategy sets. It also remains to be investigated as to whether
		learning can allow for weakening the common knowledge
		assumption.}

\bibliographystyle{IEEEtran}
\bibliography{ref}

\appendix  \label{sec:proof}

\use{{\em Proof of Lemma \ref{cexp-bound}: }}
Since $x_i^* = \Pi_{K_i} ( x_i^{*} - \gamma_i^k
		F_i(x^*;\theta^*))$, by the nonexpansivity of the Euclidean projector:
\begin{align}
	  \|x^{k+1} - x^*\|^2 \notag
& 
 \leq 
	   \sum_{i=1}^N (\|x^k_i - x_i^*\|^2 + (\gamma_i^k)^2
			\|F_i(x^k;\theta_i^k)-  F_i(x^*;\theta^*)\|^2 \notag  + (\gamma_i^k)^2
		\|w_i^k\|^2 )  \notag\\
		&   -  2\sum_{i=1}^N \gamma_i^k (x^k_i - x_i^*)^T
		(F_i(x^k;\theta_i^k)-  F_i(x^*;\theta^*))  \notag\\
  &   -  2\sum_{i=1}^N \gamma_i^k (x^k_i - x_i^*)^T w_i^k +
  2\sum_{i=1}^N (\gamma_i^k)^2 (F_i(x^k;\theta_i^k)    -  F_i(x^*;\theta^*))^T
  w_i^k.
\label{eq:grad_dim}
\end{align}
\begin{align}  \label{eq:gradient_diminishing}
\notag  &  \mbox{RHS of }\eqref{eq:grad_dim}      \leq \underbrace{\|x^k - x^*\|^2 + (\gamma_{\max}^k)^2
		\sum_{i=1}^N \|F_i(x^k;\theta_i^k)-  F_i(
				\uss{x^*};\theta^*)\|^2}_{\bf term \, 1}  + (\gamma_{\max}^k)^2 \|w^k\|^2  \\
		&\notag   \quad \underbrace{-2\sum_{i=1}^N \gamma_i^k (x^k_i - x_i^*)^T
		(F_i(x^k;\theta_i^k)-  F_i(x^k;\theta^*))}_{\bf term \, 2} \ \,  \underbrace{-2\sum_{i=1}^N \gamma_i^k (x^k_i - x_i^*)^T
			(F_i(x^k;\theta^*)-  F_i(x^*;\theta^*))}_{\bf term \, 3}\\
    &  -  2\sum_{i=1}^N \gamma_i^k (x^k_i - x_i^*)^T w_i^k   + 2\sum_{i=1}^N (\gamma_i^k)^2 (F_i(x^k;\theta_i^k)-  F_i(x^*;\theta^*))^T w_i^k.
\end{align}
By  (A\ref{assump:convex}) \usr{and using inequality $\|a+b\|^2 \leq
2\|a\|^2 + 2 \|b\|^2$}, term 1 in
(\ref{eq:gradient_diminishing}) may be bounded by leveraging the
Lipschitz continuity of $F(x;\theta)$:
\begin{align}   \label{eq:gradient_diminishing_term1}
   \notag & \usr{{\bf term \, 1} }   \leq \|x^k - x^*\|^2 + 2(\gamma_{\max}^k)^2 \sum_{i=1}^N \|F_i(x^k;\theta_i^k)-  F_i(x^k;\theta^*)\|^2  + 2(\gamma_{\max}^k)^2  \sum_{i=1}^N \|F_i(x^k;\theta^*)-  F_i(x^*;\theta^*)\|^2 \\
   \notag & \leq \|x^k - x^*\|^2 + 2(\gamma_{\max}^k)^2 \sum_{i=1}^N \|F(x^k;\theta_i^k)-  F(x^k;\theta^*)\|^2  + 2(\gamma_{\max}^k)^2   \|F(x^k;\theta^*)-  F(x^*;\theta^*)\|^2 \\
      &  \leq (1+2(\gamma_{\max}^k)^2 L_x^2) \|x^k - x^*\|^2   + 2(\gamma_{\max}^k)^2 L_{\theta}^2  \sum_{i=1}^N \|\theta_i^k-\theta^*\|^2.
\end{align}
By  (A\ref{assump:convex}), term 2 in (\ref{eq:gradient_diminishing})
	can  be bounded    by the
	Cauchy-Schwarz inequality, H\"{o}lder's inequality and the
Lipschitz continuity of $F(x;\theta)$:
\begin{align}    \label{eq:gradient_diminishing_term3}
& \quad  \notag \  \usr{{\bf term \, 2} } \\
&  \leq \notag -2\sum_{i=1}^N \gamma_i^k (x^k_i - x_i^*)^T
		(F_i(x^k;\theta_i^k)-  F_i(x^k;\theta^*))   \leq 2\gamma_{\max}^k \sum_{i=1}^N \|x_i^k - x_i^*\|
   \|F_i({x^k};\theta_i^k)-  F_i(x^k;\theta^*)\|\\
     \notag & \leq 2\gamma_{\max}^k  \|x^k - x^*\|   \sqrt{ \sum_{i=1}^N \|F_i(x^k;\theta_i^k) -
	F_i(x^k;\theta^*)\|^2}
    \\
    & \notag \leq 2\gamma_{\max}^k  \|x^k - x^*\|   \sqrt{ \sum_{i=1}^N \|F(x^k;\theta_i^k) -
	F(x^k;\theta^*)\|^2}\\
&   \leq 2 \gamma_{\max}^k  L_{\theta}  \|x^k
	- x^*\| \sqrt{\sum_{i=1}^N \|\theta_i^k-\theta^*\|^2} \leq
	\gamma_{\max}^k \mu_x \|x^k - x^*\|^2 + \gamma_{\max}^k
	\frac{L_{\theta}^2}{\mu_x} \sum_{i=1}^N \|\theta_i^k-\theta^*\|^2,
\end{align}
\uss{where the last inequality follows from the fact that $2ab\leq a^2 + b^2$.}
Term 3 in (\ref{eq:gradient_diminishing})
	can be bounded by the Cauchy-Schwarz inequality and $\gamma_i^k \leq \gamma_{\max}^k$ for all $i$:
$ \usr{{\bf term \, 3} \leq} -2\sum_{i=1}^N \gamma_i^k (x^k_i - x_i^*)^T
			(F_i(x^k;\theta^*)-  F_i(x^*;\theta^*))
  = -2\sum_{i=1}^N \gamma_{\max}^k (x^k_i - x_i^*)^T (F_i(x^k;\theta^*) -  F_i(x^*;\theta^*))
  -2\sum_{i=1}^N (\gamma_i^k-\gamma_{\max}^k) (x^k_i - x_i^*)^T
 (F_i(x^k;\theta^*) -  F_i(x^*;\theta^*))  
   \leq -2\gamma_{\max}^k\sum_{i=1}^N  (x^k_i - x_i^*)^T
	(F_i(x^k;\theta^*) -  F_i(x^*;\theta^*))  + 2	(\gamma_{\max}^k-\gamma_{\min}^k) \sum_{i=1}^N \|x_i^k - x_i^*\|
	\|F_i(x^k;\theta^*) -  F_i(x^*;\theta^*)\|.$ 
Proceeding further, we may leverage H\"{o}lder's inequality, the
Lipschitz continuity of $F(x;\theta)$ and (A\ref{assump:convex}) to obtain the following:
\begin{align}
\label{eq:gradient_diminishing_term2}
\notag     & \quad  -2\gamma_{\max}^k\sum_{i=1}^N  (x^k_i - x_i^*)^T
	(F_i(x^k;\theta^*) -  F_i(x^*;\theta^*))  \\
& \notag + 2	(\gamma_{\max}^k-\gamma_{\min}^k) \sum_{i=1}^N \|x_i^k - x_i^*\|
	\|F_i(x^k;\theta^*) -  F_i(x^*;\theta^*)\| \\
  \notag  &   \leq -2 \gamma_{\max}^k (x^k - x^*)^T (F(x^k;\theta^*) -  F(x^*;\theta^*))   + 2 (\gamma_{\max}^k-\gamma_{\min}^k)  \|x^k - x^*\| \|F(x^k;\theta^*) -  F(x^*;\theta^*)\| \\
     &   \leq -2 \gamma_{\max}^k \mu_x \|x^k - x^*\|^2        +  2	(\gamma_{\max}^k-\gamma_{\min}^k)  L_x \|x^k - x^*\|^2.
\end{align}
Combining (\ref{eq:grad_dim}) with (\ref{eq:gradient_diminishing}), (\ref{eq:gradient_diminishing_term1}), (\ref{eq:gradient_diminishing_term3}), and (\ref{eq:gradient_diminishing_term2}), we obtain
\begin{align*}
   \|x^{k+1} - x^*\|^2 & \leq (1+ 2(\gamma_{\max}^k)^2 L_x^2)\|x^k - x^*\|^2  - \gamma_{\max}^k\mu_x\|x^k - x^*\|^2 + 2 (\gamma_{\max}^k-\gamma_{\min}^k)  L_x\|x^k - x^*\|^2 \\
   &   + 2(\gamma_{\max}^k)^2 L_{\theta}^2\sum_{i=1}^N \|\theta_i^k-\theta^*\|^2 + \gamma_{\max}^k  L_{\theta}^2/\mu_x \sum_{i=1}^N \|\theta_i^k-\theta^*\|^2
    + (\gamma_{\max}^k)^2 \|w^k\|^2 \\
    &   -  2\sum_{i=1}^N \gamma_i^k (x^k_i - x_i^*)^T w_i^k +	2\sum_{i=1}^N (\gamma_i^k)^2 (F_i(x^k;\theta_i^k)-
			F_i(x^*;\theta^*))^T w_i^k\\
& =  \left(1- \gamma_{\max}^k\mu_x + 2 (\gamma_{\max}^k-\gamma_{\min}^k)  L_x +
			2(\gamma_{\max}^k)^2 L_x^2\right)  \|x^k - x^*\|^2 \\
	&   +( 2(\gamma_{\max}^k)^2 L_{\theta}^2
    + \gamma_{\max}^k  L_{\theta}^2/\mu_x ) \sum_{i=1}^N \|\theta_i^k-\theta^*\|^2  + (\gamma_{\max}^k)^2 \|w^k\|^2 \\
    & -  2\sum_{i=1}^N \gamma_i^k (x^k_i - x_i^*)^T w_i^k  + 2\sum_{i=1}^N (\gamma_i^k)^2 (F_i(x^k;\theta_i^k)-  F_i(x^*;\theta^*))^T w_i^k.
\end{align*}
{By taking conditional expectations and by recalling that $\mathbb{E}[\usr{w_i^k}
\mid {\mathcal F}_k] = 0$ and $\mathbb{E}[ \|\usr{w_i^k}\|^2 \mid {\mathcal
	F}_k] \leq \nu_x^2$, we obtain that
$\mathbb{E}\left[\|x^{k+1}-x^*\|^2\mid \mathcal{F}_k\right]  \leq
\zeta_k  \|x^k-x^*\|^2
	 +\beta_k$,
where $\zeta_k  =  1- \gamma_{\max}^k\mu_x + 2 (\gamma_{\max}^k-\gamma_{\min}^k)  L_x +
			2(\gamma_{\max}^k)^2 L_x^2$ and $\beta_k =   ( 2(\gamma_{\max}^k)^2 L_{\theta}^2
    + \gamma_{\max}^k  L_{\theta}^2/\mu_x ) \sum_{i=1}^N \|\theta_i^k-\theta^*\|^2
	    + \usr{N} (\gamma_{\max}^k)^2 \nu_x^2$.} \qed

\use{{\em Proof of Theorem \ref{thm:gradient_strongly_convergence}: }}
{From Lemma~\ref{cexp-bound}, the following holds for every $k$}:
\begin{align}\begin{aligned} \label{exp-bound-thm}
 \mathbb{E}\left[\|x^{k+1}-x^*\|^2 \mid \mathcal{F}_k\right]
& \leq
\underbrace{\left( 1- \gamma_{\max}^k\mu_x + 2 (\gamma_{\max}^k-\gamma_{\min}^k)  L_x +
			2(\gamma_{\max}^k)^2 L_x^2\right)}_{\triangleq \,
	{\zeta_k}}   \|x^k-x^*\|^2 \\
	&   + \underbrace{( 2(\gamma_{\max}^k)^2 L_{\theta}^2
    + \gamma_{\max}^k  L_{\theta}^2/\mu_x ) \sum_{i=1}^N \|\theta_i^k-\theta^*\|^2
	    + \usr{N}(\gamma_{\max}^k)^2 \nu_x^2}_{\triangleq \, \beta_k}. \end{aligned}\end{align}
{By invoking the fixed-point property given by  $\theta^* = \Pi_{\Theta} ( \theta^{*} - \alpha_{i}^k
		\nabla_{\theta} g(\theta^*))$ (see \cite{Pang03I}) and the
non-expansivity of the Euclidean projector,  we may
derive the following bound on $\|\theta_i^{k+1}-\theta^*\|^2$}:
\begin{align*} 
& \quad\ \notag   \|\theta_i^{k+1}-\theta^*\|^2   \leq \| \theta_i^{k}-\theta_i^{*} - \alpha_i^k(\nabla_{\theta} g(\theta_i^{k}) - \nabla_{\theta} g(\theta^*)) - \alpha_i^kv_i^k \|^2 \\
  & = \| \theta_i^{k}-\theta_i^{*} - \alpha_i^k(\nabla_{\theta} g(\theta_i^{k}) - \nabla_{\theta} g(\theta^*))\|^2 + (\alpha_i^k)^2 \|v_i^k \|^2  - 2 \uss{\alpha_i^k}( \theta_i^{k}-\theta_i^{*} - \alpha_i^k(\nabla_{\theta} g(\theta_i^{k}) - \nabla_{\theta} g(\theta^*)) )^T v_i^k.
\end{align*}
By taking conditional expectations, recalling that $\mathbb{E}[v_i^k
\mid \mathcal{F}_k] = 0$ and using Lemma \ref{lem:strongly_Lipschitz}, we obtain
 \begin{align}  \label{eq:stoch_strongly_optim_exp_final_theta}
 \begin{aligned}
 \mathbb{E}[\|\theta_i^{k+1}-\theta^{*}\|^2 \mid \mathcal{F}_k] &  \leq \| \theta_i^{k}-\theta^{*} - \alpha_i^k(\nabla_{\theta} g(\theta_i^{k}) - \nabla_{\theta} g(\theta^*))\|^2  + ( \alpha_i^k)^2 \mathbb{E}[\| v_i^k \|^2\mid\mathcal{F}_k] \\
   & \leq (1-2\alpha_i^k\mu_{\theta} +
	(\alpha_i^k)^2 C_{\theta}^2) \|\theta_i^{k}-\theta^{*}\|^2 + ( \alpha_i^k)^2 \nu_{\theta}^2 \\
  & \leq (1-2\alpha_{\min}^k\mu_{\theta} +
	(\alpha_{\max}^k)^2 C_{\theta}^2) \|\theta_i^{k}-\theta^{*}\|^2 + ( \alpha_{\max}^k)^2 \nu_{\theta}^2.
    \end{aligned}
\end{align}
Next, by adding $\eqref{exp-bound-thm}$ and
\eqref{eq:stoch_strongly_optim_exp_final_theta} and by invoking
(A\ref{assump:steplength_strongly}), we obtain the following:
\begin{align*} 
  \notag&\quad\  \mathbb{E}\left[\|x^{k+1}-x^{*}\|^2 \mid \mathcal{F}_k\right]  +
  \mathbb{E}\left[  \sum_{i=1}^N \|\theta^{k+1}_i-\theta^{*}\|^2\mid \mathcal{F}_k\right] \\
	 \notag  &\leq \left( 1- \gamma_{\max}^k\mu_x + 2 (\gamma_{\max}^k-\gamma_{\min}^k)  L_x +
			2(\gamma_{\max}^k)^2 L_x^2\right) \|x^{k}-x^{*}\|^2
\end{align*}
\begin{align*}
& + (1-2\alpha_{\min}^k\mu_{\theta} +
	(\alpha_{\max}^k)^2 C_{\theta}^2 + 2(\gamma_{\max}^k)^2 L_{\theta}^2
    + \gamma_{\max}^k  L_{\theta}^2/\mu_x  )  \sum_{i=1}^N\|\theta^{k}_i-\theta^*\|^2 + \usr{N} ( \gamma_{\max}^k)^2 \nu_{x}^2  +   N(\alpha_{\max}^k)^2 \nu_{\theta}^2\notag \\
\notag & \leq (1- \gamma_{\max}^k\mu_x + 2 (\gamma_{\max}^k-\gamma_{\min}^k)  L_x +
			2(\alpha_{\max}^k)^2 L_x^2\mu_x^2\mu_{\theta}^2/L_{\theta}^4  )  \|x^{k}-x^{*}\|^2 \\
  & + (1-\gamma_{\max}^k  L_{\theta}^2/\mu_x +
	(\alpha_{\max}^k)^2 C_{\theta}^2 + 2(\alpha_{\max}^k)^2  \mu_x^2\mu_{\theta}^2/L_{\theta}^2   )  \sum_{i=1}^N\|\theta^{k}_i-\theta^*\|^2 + \usr{N}( \gamma_{\max}^k)^2 \nu_{x}^2  +   N(\alpha_{\max}^k)^2 \nu_{\theta}^2\\
     \notag & \leq \left(1-\upsilon_k \gamma_{\max}^k + \beta
			 (\alpha_{\max}^k)^2 \right) \left( \|x^{k}-x^{*}\|^2 +
				 \sum_{i=1}^N\|\theta^{k}_i-\theta^*\|^2 \right)   + \delta_k,
\end{align*}
where {the second inequality results from invoking
	A\ref{assump:steplength_strongly}(c) through which $-\mu_{\theta} \alpha^k_{\rm
		min} \leq -\gamma^k_{\max} L_{\theta}^2/\mu_x$} and
		$\upsilon_k=\min\{\mu_x -
		2(\gamma_{\max}^k-\gamma_{\min}^k)L_x/\gamma_{\max}^k,
		L_{\theta}^2/\mu_{x}\}$, $\beta=\max\{2
		L_x^2\mu_x^2\mu_{\theta}^2/L_{\theta}^4  ,
		C_{\theta}^2+2\mu_x^2\mu_{\theta}^2/L_{\theta}^2 \}$, and
		$\delta_k=\usr{N}( \gamma_{\max}^k)^2 \nu_{x}^2  +
		N(\alpha_{\max}^k)^2 \nu_{\theta}^2$. {To show the
			non-summability of $(\upsilon_k \gamma_{\max}^k - \beta (
				\uss{\alpha}_{\max}^k)^2)$, we consider two cases: (i) If $\mu_x \leq
			L_{\theta}^2/\mu_x$ then $\upsilon_k = \mu_x -
		2(\gamma_{\max}^k-\gamma_{\min}^k)L_x/\gamma_{\max}^k$ and for
        \usr{some $K$, $\upsilon_k \geq \mu_x/2$ when $k > K$.}
		Consequently, $\sum_{k > K} \upsilon_k
		\gamma_{\max}^k \geq \sum_{k > K} \usr{\mu_x	\gamma_{\max}^k /2} = \infty;$ (ii) Alternately, if $\mu_x >
			L_{\theta}^2/\mu_x$, then for \usr{some $K$,} $\upsilon_k =
			L_{\theta}^2/\mu_x$ \usr{when $k > K$} and $\sum_{k > K} \upsilon_k
		\gamma_{\max}^k = \sum_{k > K} L_{\theta}^2/\mu_x
		\gamma_{\max}^k = \infty.$ Since $\uss{\alpha}_{\max}^k$ is square
		summable from (A\ref{assump:steplength_strongly}), we conclude
		that $\sum_{k=0}^{\infty} (\upsilon_k \gamma_{\max}^k - \beta (
				\uss{\alpha}_{\max}^k)^2) = \infty.$} {In addition, we
			have that
\begin{align*}&\lim_{k\to \infty} \frac{\delta_k}{{\upsilon_k} \gamma_{\max}^k -
	\beta  (\uss{\alpha}_{\max}^k)^2} = \lim_{k\to \infty} \frac{\usr{N}(\gamma_{\max}^k)^2 \nu_{x}^2  +
		N(\alpha_{\max}^k)^2 \nu_{\theta}^2 }{{\upsilon_k} \gamma_{\max}^k - \beta  (\uss{\alpha}_{\max}^k)^2} = \lim_{k\to \infty} \frac{\usr{N}( \gamma_{\max}^k) \nu_{x}^2  +
		N\frac{(\alpha_{\max}^k)^2}{\gamma_{\max}^k} \nu_{\theta}^2
		}{{\upsilon_k}  - \beta  \uss{(\alpha_{\max}^k)^2/\gamma_{\max}^k}} = 0,
\end{align*}
where the last equality results from noting that $\lim_{k \to \infty} \gamma_{\max}^k
= 0$, $\lim_{k\to \infty} (\alpha_{\max}^k)^2/\gamma_{\max}^k = 0$ and
$\lim_{k \to \infty}
\upsilon_k = c > 0.$}
Then, by invoking the super-martingale convergence theorem (Lemma
		\ref{lem:supermartingale}), we have that
$\|x^{k}-x^{*}\|^2+\sum_{i=1}^N\|\theta_i^k-\theta^*\|^2\to 0$ $a.s.$ as $k\to\infty$, which implies that
$x^k \to x^*$  and $\theta_i^k \to \theta^*$ $a.s.$ as $k\to\infty$ for all $i$. \qed

\use{{\em Proof of Theorem \ref{thm:stoch_strongly_optim_error}: }}
\usd{Throughout this proof, $\lambda_{x,\min} \triangleq {\displaystyle \min_{1\leq i \leq
N}}\{\lambda_{x,i} \}$, $\lambda_{x,\max} \triangleq {\displaystyle\max_{1\leq i \leq
		N}}\{\lambda_{x,i} \}$, $\lambda_{\theta,\min} \triangleq
		{\displaystyle \min_{1\leq i \leq
N}}\{\lambda_{\theta,i} \}$,  $\lambda_{\theta,\max} \triangleq
{\displaystyle \max_{1\leq i \leq
		N}}\{\lambda_{\theta,i} \}$. Furthermore,
$2\mu_{\theta}\lambda_{\theta,\min} > 1$,
	$\mu_x\lambda_{x,\max}-2(\lambda_{x,\max}-\lambda_{x,\min})L_x >
	1$},  $A_k \triangleq \frac{1}{2} \|x^k - x^*\|^2$, and $a_k \triangleq
\mathbb{E}[A_k]$. Then, $A_{k+1}$ may be bounded as follows by using the
non-expansivity of the Euclidean projector:
\begin{align}    \label{alg:stoch_grad_grad_strongly_constant_A}
   A_{k+1} 
   &\notag = \frac{1}{2} \sum_{i=1}^N \|x_i^{k+1} - x_i^*\|^2 = \frac{1}{2} \sum_{i=1}^N \|\Pi_{K_i} (x_i^k - \gamma_{i}^k \nabla_{x_i}
				f_i(x^k;{\theta_i^k},\xi_i^k)) - \Pi_{K_i} (x_i^*)\|^2  \\
 & \notag \leq \frac{1}{2} \sum_{i=1}^N \| x_i^k - \gamma_{i}^k \nabla_{x_i}
				f_i(x^k;{\theta_i^k},\xi_i^k) - x_i^* \|^2   \leq\frac{1}{2}\sum_{i=1}^N \| x_i^{k}-x_i^{*} -
 \gamma_i^k(F_i(x^{k};\theta_i^{k}) + w_i^k)) \|^2 \\
	 &   = A_k + \frac{1}{2} \sum_{i=1}^N(\gamma_i^k)^2 \|F_i(x^{k};\theta_i^{k}) + w_i^k \|^2  - \sum_{i=1}^N \gamma_i^k (x_i^{k}-x_i^{*})^T ( F_i(x^{k};\theta_i^{k}) + w_i^k ).
\end{align}
Note that $\mathbb{E}[(x_i^{k}-x_i^{*})^Tw_i^k] \usr{= \mathbb{E}[\mathbb{E}[ (w_i^k)^T(x_i^k-x_i^*)\mid \Fscr_k]]} = 0$.
By taking expectations on both sides of
\usd{\eqref{alg:stoch_grad_grad_strongly_constant_A}} and by invoking the
bounds
$\mathbb{E}[\|F_i(x^{k};\theta_i^{k}) + w_i^k \|^2]\leq M^2/N$  and
$\mathbb{E}[\|\nabla_{\theta} g(\theta_i^{k})+ v_i^k \|^2]\leq M_{\theta}^2$,
it follows that
\begin{align} \label{alg:stoch_grad_grad_strongly_constant_a}
  a_{k+1}
   & \leq a_k + \frac{1}{2} (\gamma_{\max}^k)^2 M^2 - \sum_{i=1}^N \gamma_i^k \mathbb{E}[(x_i^{k}-x_i^{*})^T  F_i(x^{k};\theta_i^{k})].
\end{align}
By \eqref{eq:gradient_diminishing_term3} and \eqref{eq:gradient_diminishing_term2}, the last term (including the negative sign) in \eqref{alg:stoch_grad_grad_strongly_constant_a} can be bounded by
\begin{align}
 \label{alg:stoch_grad_grad_strongly_constant_a_term2}
 \notag & \quad\ -\sum_{i=1}^N \gamma_i^k \mathbb{E}[(x_i^{k}-x_i^{*})^T ( F_i(x^{k};\theta_i^{k}) -  F_i(x^{k};\theta^{*})) ]    \\
\notag &\quad\   - \sum_{i=1}^N \gamma_i^k \mathbb{E}[(x_i^{k}-x_i^{*})^T ( F_i(x^{k};\theta^{*}) -  F_i(x^{*};\theta^{*})) ]    - \sum_{i=1}^N \gamma_i^k \mathbb{E}[(x_i^{k}-x_i^{*})^T F_i(x^{*};\theta^{*}) ]   \\
\notag & \leq -\sum_{i=1}^N \gamma_i^k \mathbb{E}[(x_i^{k}-x_i^{*})^T ( F_i(x^{k};\theta_i^{k}) -  F_i(x^{k};\theta^{*})) ]   - \sum_{i=1}^N \gamma_i^k \mathbb{E}[(x_i^{k}-x_i^{*})^T ( F_i(x^{k};\theta^{*}) -  F_i(x^{*};\theta^{*})) ] \\
   & \leq   \gamma_{\max}^k \mu_x a_k + \gamma_{\max}^k  L_{\theta}^2/(2\mu_x)  \sum_{i=1}^N \mathbb{E}[\|\theta_i^k-\theta^*\|^2]  -    2\gamma_{\max}^k \mu_x  a_k     +
	2(\gamma_{\max}^k-\gamma_{\min}^k)  L_x a_k.
\end{align}
Combining \eqref{alg:stoch_grad_grad_strongly_constant_a} and \eqref{alg:stoch_grad_grad_strongly_constant_a_term2}, we get
\begin{align*}
   a_{k+1}
    &\leq(1- \gamma_{\max}^k \mu_x        +
	2(\gamma_{\max}^k-\gamma_{\min}^k)  L_x )a_k   + \frac{1}{2} (\gamma_{\max}^k)^2 M^2     + \gamma_{\max}^k  L_{\theta}^2/(2\mu_x)  \sum_{i=1}^N \mathbb{E}[\|\theta_i^k-\theta^*\|^2].
\end{align*}
Suppose $\alpha_i^k=\lambda_{\theta,i}/k$ for all $i$. Since the function $g(\theta)$ is strongly convex, we can use the standard rate estimate (cf. inequality (5.292) in \cite{Shapiro09lecturesSA}) to get the following
\begin{align} \label{eq:error_theta_Q}
  \mathbb{E}[\|\theta_i^{k} - \theta^{*} \|^2 ] \leq
  \frac{\usd{Q_{\theta}}}{k},
\end{align}
where $Q_{\theta}\triangleq \max \{ \lambda_{\theta,\max}^2 M_{\theta}^2 (2\mu_{\theta}\lambda_{\theta,\min} - 1)^{-1}$, $\max_i\mathbb{E}[\|\theta_i^{0} - \theta^{*} \|^2 ] \}$ with $\lambda_{\theta,\min}>1/(2\mu_{\theta})$.
Suppose $\gamma_i^k=\lambda_{x,i}/k$, allowing us to claim the
following:
\begin{align*}
     a_{k+1}
      &\leq\left(1-\frac{\mu_x\lambda_{x,\max}-2(\lambda_{x,\max}-\lambda_{x,\min})L_x}{k}\right)a_k + \frac{\lambda_{x,\max}^2}{2k^2}
	  \left(M^2 + \frac{L_{\theta}^2 N \usd{Q_{\theta}}}{\lambda_{x,\max}\mu_x
		  }\right),
\end{align*}
By assuming that $\mu_x\lambda_{x,\max}-2(\lambda_{x,\max}-\lambda_{x,\min})L_x>1$, the result follows by
observing that
{$\mathbb{E}[\|x^{k} - x^{*} \|^2 ] \leq
   \frac{\usd{Q_{x,{\theta}}}}{k},$}
where
\usd{$Q_{x,{\theta}}\triangleq
	\max\left\{ \frac{\lambda_{x,\max}^2M^2 + \lambda_{x,\max}^2
	  L^2_{\theta}NQ_{\theta}}{(\mu_x\lambda_{x,\max}-2(\lambda_{x,\max}-\lambda_{x,\min})L_x-1)}, \mathbb{E}[\|x^{0} -
x^{*} \|^2 ]\right\}.$} \qed

{\em Proof of Proposition \ref{lem:best_response_rand}: }
It suffices to show that given $p(X^k;\theta^*,\xi^k)$ and
$\{\vthetabar_i^k\}_{i=1}^N$, the variational inequality VI$({\cal Z}, F^{k+1} \usr{+ \epsilon^k\mathbf{I}} )$ has a unique solution for each $i$.
Now, for simplicity, we ignore the superscript $k$  for all variables.
Given $p$, $\vthetabar_i$, $i$ and $k$,  let $H(z_i)$ denote the Jacobian
matrix $\nabla F(z_i)$ of $F$ at $z_i\in \cal Z$.
We will proceed to show that $H(z_i)$ is a $\mathbf{P}$-matrix for all
$z_i\in \cal \widetilde{Z}$ in part (a) and a $\mathbf{P_0}$-matrix for all
$z_i\in \cal \widetilde{Z}$ in part (b) where $\cal Z \subset \cal \widetilde{Z}$ and
$\cal \widetilde{Z}$ is a rectangle.  Then, by invoking Proposition 3.5.9 in \cite{Pang03I}, the associated
mapping $F$ is $\textbf{P}$-mapping on $\cal \widetilde{Z}$ in part (a) and a
$\mathbf{P}_0$-mapping on $\cal \widetilde{Z}$ in part (b). Consequently, by
Theorem 3.5.15 in \cite{Pang03I}, the regularized variational inequality VI$({\cal Z}, F^{k+1} \usr{+ \epsilon^k\mathbf{I}})$  has a unique solution in both parts (a)
	and (b). Specifically, we employ a rectangular $\widetilde{\cal Z}$ defined as
	$  \widetilde{\cal Z} \triangleq [0,\infty)^N \times \Theta,$
	where $\Theta$ is a compact set in $(0,\infty)$.
\noindent {\bf (a)}  Given $z_i\in \widetilde{\cal Z}$, let $H_i$ denote $H(z_i)$. Then,
$H_i = \left(
  \begin{array}{cc}
    A_i & B \\
    C & D
  \end{array}
\right)$,
where
$A_i=b^*(I+ee^T)+E_i$, $B=-\frac{1}{k+1}e$, $C = - b^* e^T$, $D=\frac{1}{k+1}$,
$e$ denotes the
column of ones in $\Real^N$,
$E_i$ is an $N\times
N$ diagonal matrix with $c_j''(x_{ij})$ as its $j$th diagonal entry.
Since, the nonnegativity of $c_j''(x_{ij})$ follows from the convexity
of costs, $E_i$ is a nonnegative diagonal matrix and is therefore
positive semidefinite.  Recall that the sum of a diagonal positive semidefinite matrix and a
$\mathbf{P}$-matrix is a $\mathbf{P}$-matrix and it suffices
to show that $H_i$ is a $\mathbf{P}$-matrix when $E_i=0$. This amounts to
showing that the principal minors of $H$ are positive.

Since $A_i$ and $D$ are $\mathbf{P}$-matrices, we only consider the
principal submatrix $H_{\alpha}$ of $H_i$,
where $\alpha \subseteq
\{1,\ldots,N \}$ is a nonempty index set and $H_{\alpha}$ is given by
$H_{\alpha}
 \triangleq
\left(
  \begin{array}{cccc}
    A_{\alpha} &  B_{\alpha} \\
    C_{\alpha}  & D
  \end{array}
\right),$
where
     $A_{\alpha}  = b^*(I_{n_{\alpha}}+e^{n_{\alpha}}(e^{n_{\alpha}})^T),
     B_{\alpha}  = -\frac{1}{k+1}e^{n_{\alpha}},
     C_{\alpha}  = - b^*  (e^{n_{\alpha}})^T$,
and
$I_{n_{\alpha}}$ and $e^{n_{\alpha}}$ denote the identity matrix and the
column of ones in $\Real^{n_{\alpha}\times n_{\alpha}}$ and $\Real^{n_{\alpha}}$, respectively,
with $n_{\alpha}=|\alpha|$.
Since
$A_{\alpha}^{-1}=\frac{1}{b^*}\left(I_{n_{\alpha}}-\frac{1}{n_{\alpha}+1}e^{n_{\alpha}}(e^{n_{\alpha}})^T\right)$, we have
  $C_{\alpha}A_{\alpha}^{-1}B_{\alpha}
     = \frac{1}{k+1}(e^{n_{\alpha}})^T
	\left(I_{n_{\alpha}}-\frac{1}{n_{\alpha}+1}e^{n_{\alpha}}(e^{n_{\alpha}})^T\right) e^{n_{\alpha}}$ $= \frac{1}{k+1}\left(n_{\alpha} -
			\frac{n_{\alpha}^2}{n_{\alpha}+1}\right) =
	\frac{1}{k+1}\left(\frac{n_{\alpha}}{n_{\alpha}+1}\right).$
It follows that
$D-C_{\alpha}A_{\alpha}^{-1}B_{\alpha}  = \frac{1}{k+1}-\frac{1}{k+1}
\left(\frac{n_{\alpha}}{n_{\alpha}+1}\right) =
\frac{1}{k+1}\left(\frac{1}{n_{\alpha}+1}\right)>0.$
Since $\det(A_{\alpha})>0$,
we have
$\det(H_{\alpha})=\det(A_{\alpha})\det\left(D-C_{\alpha}A_{\alpha}^{-1}B_{\alpha}\right)>0$ for all
$\alpha \subseteq\{1,\ldots,N \}$ with $\alpha\neq\emptyset$. Therefore, $H$ is a $\mathbf{P}$-matrix.

\noindent {\bf (b)}  Analogous to our approach for (a), we consider a matrix $H_i$,
	given by $H_i= \nabla F(z_i)$. Then,
$H_i =
\left(
  \begin{array}{cc}
    A_i & B_i \\
    C_i & D_i
  \end{array}
\right)$,
where
$A_i=\bhat_i(I+ee^T)+E_i$, $B_i=\frac{1}{k+1}(\usd{\bf x}_i+(e^T \usd{\bf x}_i)e)$, $C_i = \bhat_i e^T$, and $D_i=\frac{1}{k+1}(e^T \usd{\bf x}_i)$,
where
$\bhat_i=\frac{1}{k+1}b_i+\frac{k}{k+1} \bbar_i$,
$\usd{\bf x}_i=(x_{i1},\ldots,x_{iN})^T$, $e$ denotes the
column of ones in $\Real^N$,
and $E_i$ is an $N\times
N$ diagonal matrix with $c_j''(x_{ij})$ as its $j$th diagonal entry.
Recall that the sum of a diagonal positive semidefinite matrix and a
$\mathbf{P}_0$-matrix is a $\mathbf{P}_0$-matrix. As in (a), it suffices
to show that $H$ is a $\mathbf{P}_0$-matrix when $E_i=0$.

Since $A_i$ and $D_i$ are $\mathbf{P}_0$-matrices, we restrict our
attention to the principal submatrix $H_{\alpha}$ of $H_i$,
where $\alpha \subseteq
\{1,\ldots,N \}$ is a nonempty index set, and $H_{\alpha}$ is given by
$H_{\alpha}
\triangleq
\left(
  \begin{array}{cccc}
    A_{\alpha} &  B_{\alpha} \\
    C_{\alpha}  & D_i
  \end{array}
\right),$
where
$     A_{\alpha}   = \bhat_i(I_{n_{\alpha}}+e^{n_{\alpha}}(e^{n_{\alpha}})^T),
     B_{\alpha}  = \frac{1}{k+1}(\usd{\bf x}_{\alpha}+(e^T x_i)e^{n_{\alpha}}),
     C_{\alpha}  = \bhat_i (e^{n_{\alpha}})^T,$
and
$I_{n_{\alpha}}$ and $e^{n_{\alpha}}$ denote the identity matrix and the
column of ones in $\Real^{n_{\alpha}\times n_{\alpha}}$ and $\Real^{n_{\alpha}}$, respectively,
with $n_{\alpha}=|\alpha|$. Then, the following hold:
(1) If $\bhat_i=0$, then $A_{\alpha}=0$ and $C_{\alpha}=0$, which implies $\det(H_{\alpha})=0$;
(2) If $\bhat_i>0$, then $A_{\alpha}^{-1}=\frac{1}{\bhat_i}(I_{n_{\alpha}}-\frac{1}{n_{\alpha}+1}e^{n_{\alpha}}(e^{n_{\alpha}})^T)$.
So, we have
\begin{align*}
 &  \quad\,  C_{\alpha}A_{\alpha}^{-1}B_{\alpha}     = \frac{1}{k+1}(e^{n_{\alpha}})^T
	\left(I_{n_{\alpha}}-\frac{1}{n_{\alpha}+1}e^{n_{\alpha}}(e^{n_{\alpha}})^T\right)  (x_{\alpha}+(e^T x_i)e^{n_{\alpha}}) \\
    & = \frac{1}{k+1}((e^{n_{\alpha}})^T -
			\frac{n_{\alpha}}{n_{\alpha}+1}\left(e^{n_{\alpha}})^T \right) (x_{\alpha}+(e^T x_i)e^{n_{\alpha}})   = \frac{1}{(k+1)(n_{\alpha}+1)} \left((e^{n_{\alpha}})^T
			x_{\alpha} + n_{\alpha}(e^T x_i)\right).
\end{align*}
\begin{align*}
\implies   D_i-C_{\alpha}A_{\alpha}^{-1}B_{\alpha} & = \frac{1}{k+1} e^T x_i -
\frac{1}{(k+1)(n_{\alpha}+1)} \left((e^{n_{\alpha}})^T x_{\alpha} +
		n_{\alpha}(e^T x_i)\right) \\
&  = \frac{1}{(k+1)(n_{\alpha}+1)} \left(e^T x_i - (e^{n_{\alpha}})^T
		x_{\alpha}\right) \geq 0.
\end{align*}
Since $\det(A_{\alpha})>0$,
we have
$\det(H_{\alpha})=\det(A_{\alpha})\det\left(D_i-C_{\alpha}A_{\alpha}^{-1}B_{\alpha}\right)\geq 0$.

Therefore, $\det(H_{\alpha}) \geq 0$ for all nonempty
$\alpha \subseteq\{1,\ldots,N \}$,
implying that $H_i$ is a $\mathbf{P_0}$-matrix.\qed

{\em Proof of Lemma \ref{mon-Lip}:}
Let $g(x)=(c'_1(x_1),\dots,c'_N(x_N))^T$ and $e=(1,\ldots,1)^T$. Then, we have
$F(x)=g(x)+b^*(x+Xe)-a^*e$, where $X=\sum_{i=1}^N x_i$.
Note that $g(x)$ is monotone in $x$. Thus, we have for $x,y\in K$
\begin{align*}
   (F(x)-F(y))^T(x-y)& = (g(x)-g(y))^T (x-y) + b^* (x-y)^T(x-y) + b^* (X-Y)e^T(x-y) \\
   &\geq b^* (x-y)^T(x-y) + b^* (X-Y)^T(X-Y)
   \geq b^* \|x-y\|^2.
\end{align*}
This implies that  $F(x)$ is strongly monotone in $x$ with constant
		$b^*.$
Note that $g(x)$ is Lipschitz continuous on $K$ with constant $M$, where $M \triangleq \max_i\{M_i\}$. The Lipschitz continuity of $F(x)$ is easily shown:
	\begin{align*}
		\|F(x)-F(y)\| & = \|g(x)-g(y)\| + b^* \|x-y\| + b^* \|(X-Y)e\|\\
& \leq M \|x-y\|  + b^* \|x-y\| + b^* \|ee^T\|\|x-y\|
 = L\|x-y\|,
	\end{align*} where $L=M+b^*+b^*\|ee^T\|$. It follows
	that $F(x)$ is Lipschitz continuous with constant $L$.\qed

\use{{\em Proof of Proposition \ref{prop:unique}:}}
From Lemma~\ref{mon-Lip}, the associated variational
	inequality VI$(K,F)$ has a strongly monotone mapping $F(x)$ over
	$K$. Consequently, VI$(K,F)$ admits a unique
	solution~\cite{Pang03I}. \qed

{\em Proof of Proposition \ref{prop:br}:}
Consider $\theta_1,\theta_2\in\Theta$ and let
$F_i(\cdot):=F(\cdot,\theta_i)$, $i = 1, 2$. Let $x_i$ be a  solution of
VI$(K,F_i)$ for $i = 1,2$. By the assumption of strong monotonicity on
			 the map, we have that
\begin{align}\label{temp111}
   (x_1-x_2)^T(F_1(x_1)-F_1(x_2)) \geq c\|x_2-x_1\|^2,
\end{align}
for some constant $c>0$ (assumed to be independent of $\theta_1$).
Since $x_1$ is a solution of VI$(K,F_1)$, it follows that $(x_2-x_1)^TF_1(x_1)\geq 0$, which together with \eqref{temp111}
implies
\begin{align}\label{temp112}
   (x_2-x_1)^TF_1(x_2) \geq c\|x_2-x_1\|^2.
\end{align}
We may express \eqref{temp112}  as
   $(x_2-x_1)^T(F_1(x_2)-F_2(x_2)+F_2(x_2)) \geq c\|x_2-x_1\|^2.$
Now since $x_2$ is the
solution of VI$(K,F_2)$, it follows that  $(x_2-x_1)^TF_2(x_2)\leq 0$.
Consequently we obtain
\begin{align}\label{temp113}
  \|x_2-x_1\|\|F_1(x_2)-F_2(x_2)\| \geq (x_2-x_1)^T(F_1(x_2)-F_2(x_2))\geq c\|x_2-x_1\|^2.
\end{align}
By Lipschitz continuity of $F(x,\theta)$ (assuming it is uniform in
		$x$), we have that $\|F_1(x_2,\theta_1)-F_2(x_2,\theta_2)\|\leq
L_{\theta}\|\theta_2-\theta_1\|$,
and hence by  \eqref{temp113}
$L_{\theta}\|x_2-x_1\|\|\theta_2-\theta_1\|\geq c\|x_2-x_1\|^2.$
 It follows that $\|x_2-x_1\|\leq L_{\theta}c^{-1}\|\theta_2-\theta_1\|$.


To show (b), let $x(\theta_i,\epsilon_j)$ be the solution of
VI$(K,G_{ij}(\cdot))$, where
$G_{ij}(\cdot)=F(\cdot,\theta_i)+\epsilon_j \mathbf{I}$.
We begin by applying the triangle inequality to obtain that
$\|x(\theta_1,\epsilon_1) -x(\theta_2,\epsilon_2)\| \leq
\|x(\theta_1,\epsilon_1)-x(\theta_2,\epsilon_1)\| +
\|x(\theta_2,\epsilon_1) - x(\theta_2,\epsilon_2)\|.$
Since $G_{i1}$ is strongly monotone in $x$ with constant $c+\epsilon_1$
and Lipschitz continuous in $\theta$ with constant $L_{\theta}$, respectively, we have
that the first term is bounded by $L_{\theta}(c+\epsilon_1)^{-1}\|\theta_2-\theta_1\|$ as a result from part (a).
 Before proceeding, the
Lipschitz continuity of $F(x;\theta)+\epsilon I$ with respect to
$\epsilon$ can be obtained as
$$ \|(F(x;\theta)+\epsilon_2 x)-(F(x;\theta)+\epsilon_1x)\|
\leq \|x\| \|\epsilon_1-\epsilon_2\| \leq D\|\epsilon_1 - \epsilon_2\|.$$
Since $G_{2j}$ is strongly monotone in $x$ with constant $c+\epsilon_j$ and Lipschitz continuous in $\epsilon$ with constant $D$, respectively, we have
that the second term is bounded by
$D(c+\epsilon_1)^{-1}\|\epsilon_2-\epsilon_1\|$ as a result from part
(a).  Consequently, we obtain that $$\|x(\theta_1,\epsilon_1) -x(\theta_2,\epsilon_2)\| \leq
L_{\theta}(c+\epsilon_1)^{-1}\|\theta_2-\theta_1\| +
D(c+\epsilon_1)^{-1}\|\epsilon_2-\epsilon_1\|.$$
The Lipschitz continuity of $x(\theta,\epsilon)$ with respect to its
parameters follows.\qed

\use{{\em Proof of Theorem \ref{thm:br_sto}:}}
Suppose $k\geq0$. At the $k$th iteration, $\tilde{p}_i^k$
is a function of $\hat{\theta}_i^{k+1}$, which is a function of
$\bar{\vartheta}_i^k$. Consequently, the \us{fixed-point} problem
\eqref{best_response_rand} is a
function of $\bar{\vartheta}_i^k$. Since \eqref{best_response_rand}
 has a unique solution (Prop.~\ref{lem:best_response_rand}), it follows that $x_{i,\cdot} =
x_{j,\cdot}$ for $i \neq j$ and
$x_{ij}^k = x_{jj}^k$. Therefore, Given $p(X^k;\theta^*,\xi^k)$ and $\{\vthetabar_i^k\}_{i=1}^N$,the solution $(\usd{\bf x}_{i}^{k+1},\theta_i^{k+1})$ to
\eqref{best_response_rand} satisfies $x_{ij}^{k+1}=x_{jj}^{k+1}$ for all
$i,j$.  Thus, for all $k \geq 0$ and all $i$, we have that
\begin{align*}
    & p(X^k;\uss{\theta^*},\xi^k)  =
    \begin{cases}
        (a^*+\xi^{k})-b^*\sum_{j=1}^N x_{jj}^{k}  = (a^*+\xi^{k})-b^*\sum_{j=1}^N x_{ij}^{k}, &\quad \mbox{under }~
				\textbf{(A\ref{assump:ab}a)}, \\
        a^*-(\theta^*+\xi^{k})\sum_{j=1}^N x_{jj}^{k}  =
		a^*-(\theta^*+\xi^{k})\sum_{j=1}^N x_{ij}^{k},  &\quad \mbox{under }~
				\textbf{(A\ref{assump:ab}b)}.
     \end{cases}
\end{align*}
Since for all $k \geq 0$ and all $i$,
\begin{align*}
    \vtheta_i^k= \begin{cases}
			 	p(X^k;\uss{\theta^*},\xi^k)+b^*X_i^k, &\quad \mbox{under }~
				\textbf{(A\ref{assump:ab}a)},  \\
				(a^*-p(X^k;\uss{\theta^*},\xi^k))/X_i^k,  & \quad \mbox{under }~
				\textbf{(A\ref{assump:ab}b)}. \end{cases}
\end{align*}
we have $\vtheta_i^k=\theta^*+\xi^{k}$ for all $i$. As a result, after $k$ \us{iterative fixed-point} steps, we obtain $k$ samples
$\{\theta^*+\xi^1, \hdots, \theta^*+\xi^{k}\}$ of the estimated parameter.
Since for all $k \geq 0$ and all $i$,
$\usd{\bar{\vartheta}_i^{k+1}=\frac{k\vthetabar_i^{k}+\vtheta_i^{k+1}}{k+1}}$,
the sample mean of the estimated parameter is given by
$\vthetabar_i^{k}$, i.e.,
		\begin{align} \label{vtheta-savg}
		\vthetabar_i^{k}= \frac{\sum_{l=1}^{k} (\theta^*+
				\xi^l)}{{k}}.\end{align}
Therefore, $\vthetabar_i^{k} \to \theta^*$ $a.s.$ as $k \to \infty$, which implies by the boundedness of $\{\theta_i^{k}\}$ that for all $i$
\begin{align*}
    \thetahat_i^{k+1} = \frac{1}{k+1}\theta_i^{k+1}+\frac{k}{k+1}
	\vthetabar_i^k \to \theta^*  \quad\, \textrm{\textit{a.s.}\,\,  as } k
	\to \infty,
\end{align*}
by the strong law of large numbers.  By Proposition~\ref{prop:br}, $\usd{\bf x}^{k+1}_i=x^{k+1}_i(\thetahat_i^{k+1},\epsilon^k)$ is a continuous function of $(\thetahat_i^{k+1},\epsilon^k)$, and $\usd{\bf x}^{k+1}_i(\theta^*,0)=x^*$.
Therefore, $\usd{\bf x}^{k+1}_i \to x^*$ $a.s.$ as $k \to \infty$. \qed

{\em Proof of Lemma \ref{monotone-map}:}
\noindent (a) Strict monotonicity of $F(x)$ is implied by the positive
definiteness of the Jacobian $\nabla F(x).$ This is given by $
	\nabla F(x) = J_1 +J_2 +J_3$, where $J_2 = 2b^*\sigma X^{\sigma - 1}
	ee^T$ and
	 $J_1 =\pmat{ c_1''(x_1) & \\ &  \ddots \\ & & c_N''(x_N)} \mbox{ and
		 }
	, 
J_3= b^*\sigma (\sigma	-1) X^{\sigma -2} \pmat{
										\frac{X}{\sigma -1} + x_1 &
											\hdots & x_1 \\
											\vdots & \ddots & \vdots \\
											x_N & \hdots &
											\frac{X}{\sigma-1} + x_N}.$
	Since $c_i(x_i)$ is a convex function in $x_i$ for all $i$, $J_1$ is a positive semidefinite matrix.  $J_2$,
	compactly stated as $2b^* \sigma X^{\sigma-1} ee^T$, is also a positive
	semidefinite matrix.  As a consequence, positive definiteness of
	$\nabla F(x)$  follows from the diagonal dominance of the following
	matrix:
\begin{align*}
	 b^*\sigma (\sigma
											-1) X^{\sigma -2} \pmat{
										\frac{X}{\sigma -1} + x_1 &
											\hdots & {1\over 2} (x_1+x_N) \\
											\vdots & \ddots & \vdots \\
											{1\over2}(x_N+x_1) & \hdots &
											\frac{X}{\sigma-1} + x_N}.
	\end{align*}
	By a minor rearrangement, it suffices to show the diagonal dominance
	of the following:
\begin{align*}
	& b^*\sigma (\sigma
											-1) X^{\sigma -2}   \pmat{
										\frac{X_{-1}}{\sigma -1} +
											(1+\frac{1}{(\sigma-1)})x_1 &
											\hdots & {1\over 2} (x_1+x_N) \\
											\vdots & \ddots & \vdots \\
											{1\over2}(x_N+x_1) & \hdots &
											\frac{X_{-N}}{\sigma-1} +
												(1+\frac{1}{(\sigma-1)})x_N},
	\end{align*}
	where $X_{-j} \triangleq \sum_{i \neq j} x_i$.  The result follows
	by noting that
	\begin{align*}
		\left(1+\frac{1}{(\sigma-1)}\right) > \frac{(N-1)}{2} \mbox{ or }
			\frac{2\sigma}{\sigma - 1} > N - 1 \mbox{ or } N < \frac{
				3\sigma - 1}{\sigma - 1},
	\end{align*}
\usr{and}
 \usr{$\frac{1}{(\sigma-1)} \geq \frac{1}{2} \mbox{ or } 1< \sigma \leq 3.$}
\noindent (b)  For $x,y\in K$,
    $(x-y)^T(F(x)-F(y)) = \int_0^1(x-y)^T \nabla F(y+\alpha(x-y))
	(x-y)d\alpha.$
Let $\tilde{x}=y+\alpha(x-y)$ an $\tilde{X}=\sum_{i=1}^N{\tilde{x}_i}$.
Akin to $\nabla F(x)$, 
$\nabla F(y+\alpha(x-y))=\tilde{J}_1+\tilde{J}_2+\tilde{J}_3$,
where $\tilde{J}_1$ and $\tilde{J}_2$ are positive semidefinite, and $\tilde{J}_3=b^*\sigma (\sigma	-1) \tilde{X}^{\sigma -2} \tilde{J}_4 $,
where
\begin{align*}
	 \tilde{J}_4 & = \pmat{
										\frac{\tilde{X}_{-1}}{\sigma -1} +
											(1+\frac{1}{(\sigma-1)})\tilde{x}_1 &
											\hdots & {1\over 2} (\tilde{x}_1+\tilde{x}_N) \\
											\vdots & \ddots & \vdots \\
											{1\over2}(\tilde{x}_N+\tilde{x}_1) & \hdots &
											\frac{\tilde{X}_{-N}}{\sigma-1} +
												(1+\frac{1}{(\sigma-1)})\tilde{x}_N}
\end{align*}
\begin{align*}
    & = \pmat{
										\frac{\tilde{X}_{-1}}{\sigma -1} +
											\frac{N-1}{2}\tilde{x}_1 &
											\hdots & {1\over 2} (\tilde{x}_1+\tilde{x}_N) \\
											\vdots & \ddots & \vdots \\
											{1\over2}(\tilde{x}_N+\tilde{x}_1) & \hdots &
											\frac{\tilde{X}_{-N}}{\sigma-1} +
												\frac{N-1}{2}\tilde{x}_N} 	+ \left(\frac{\sigma}{\sigma-1}- \frac{N-1}{2}\right) I_N  \triangleq \tilde{J}_5 +\rho I_N,
	\end{align*}
 where $\tilde{J}_5$ is a positive semidefinite matrix and
 $\rho=(1+\frac{1}{\sigma-1}- \frac{N-1}{2})>0$.  Therefore,
\begin{align*}
 &\quad      (x-y)^T(F(x)-F(y))  \geq \int_0^1(x-y)^T  \tilde{J}_3
	(x-y)d\alpha   \geq   b^*\sigma (\sigma	-1)  \int_0^1(x-y)^T
	\tilde{X}^{\sigma -2}\tilde{J}_4 (x-y)d\alpha \\
& \geq   b^*\sigma (\sigma	-1)  \eta^{\sigma -2} \int_0^1(x-y)^T  (\tilde{J}_5 +\rho I_N) (x-y)d\alpha \geq  b^*\sigma (\sigma	-1)  \eta^{\sigma -2} \rho \|x-y\|^2,
\end{align*}
implying the strong monotonicity of $F$.\qed

\use{{\em Proof of Corollary \ref{cor:br}:}}
By Proposition~\ref{prop:br}, it suffices to show that $F(x;\theta)$ is Lipschitz
		continuous in $x$ for all $\theta \in \Theta$.  For $\theta \in
		\Theta$, and $x,y\in K$, we have that
\begin{align*}
    \|F(x;\theta)-F(y;\theta)\|  & = \left\|\int_0^1 \nabla
	F(y+\alpha(x-y);\theta) (x-y) d\alpha\right\|  \quad \textrm{for some } \alpha\in(0,1) \\
    & \leq \int_0^1 \left\|\nabla F(y+\alpha(x-y);\theta)\right\|\|x-y\| d\alpha \leq \int_0^1 r\|x-y\| d\alpha = r\|x-y\|,
       \end{align*}
which implies the Lipschitz
		continuity in $x$ of the mapping $F$. \qed

\use{{\em Proof of Proposition \ref{prop:br_non}:}}
By Lemma \ref{monotone-map},  $F(x;\theta)$ is a strongly
	monotone mapping over $K$ for all $\theta \in \Theta$. By definition of $F$, $F(x;\theta)$ is Lipschitz
		continuous in $\theta$ for all $x \in K$. By definition of $\nabla F$ and boundedness of $x\in K$,
$\nabla F(x;\theta)$ is bounded for $x\in K$ and $\theta\in \Theta$. Then, the conclusion follows from Corollary \ref{cor:br}. \qed

{\em Proof of Proposition \ref{prop:br_non_unique}:}
Given $p$, $\thetabar_i$, $i$ and $k$, let $H(z_i)$ denote the Jacobian matrix $\nabla F(z_i)$ of the mapping
$F$ at $z_i\in \widetilde{\cal Z}$.
Then, as in Proposition \ref{lem:best_response_rand},
it suffices to show that $H(z_i)$ is a $\mathbf{P}$-matrix for all
$z_i\in \widetilde{\cal Z}$.
Given $z_i\in \widetilde{\cal Z}$, let $H=H(z_i)$. Then,
$H=H(z_i)=
\left(
  \begin{array}{cc}
    A_i & B \\
    C_i & D
  \end{array}
\right)$,
where
$A_i=\sigma b^* (X_i)^{\sigma-2} \left[X_i(I+ee^T)+(\sigma-1) \usd{\bf x}_i e^T \right] + E_i$, $B=-\frac{1}{k+1}e$, $C_i = - \sigma b^* (X_i)^{\sigma-1} e^T$, and $D=\frac{1}{k+1}$,
where $X_i=\sum_{j=1}^N x_{ij}$, $\usd{\bf x}_i=(x_{11},\ldots,x_{1N})^T$, and
$E_i$ is an $N\times
N$ diagonal matrix with $c_j''(x_{ij})$ as its $j$th diagonal entry.
It suffices
to show that $H$ is a $\mathbf{P}$-matrix when $E_i=0$.

If $N < \frac{ 3\sigma - 1}{\sigma - 1}$, then $A_i$ is positive
semidefinite by Lemma~\ref{monotone-map}. Therefore, we only consider the principal submatrix $H_{\alpha}$ of $H$,
where $\alpha \subseteq
\{1,\ldots,N \}$ is a nonempty index set, and $H_{\alpha} \triangleq
\left(
  \begin{array}{cccc}
    A_{\alpha} &  B_{\alpha} \\
    C_{\alpha}  & D
  \end{array}
\right),$
where
\us{$
     A_{\alpha}   = \sigma b^* (X_i)^{\sigma-1} \left[ I_{n_{\alpha}}+e^{n_{\alpha}}(e^{n_{\alpha}})^T \right] + \sigma(\sigma-1) b^* (X_i)^{\sigma-2}\usd{\bf x}_{\alpha}(e^{n_{\alpha}})^T ,
     B_{\alpha}   = -\frac{1}{k+1}e^{n_{\alpha}},
     C_{\alpha}  = - \sigma b^* (X_i)^{\sigma-1}  (e^{n_{\alpha}})^T,
$}
and
$I_{n_{\alpha}}$ and $e^{n_{\alpha}}$ denote the identity matrix and the
column of ones in $\Real^{n_{\alpha}\times n_{\alpha}}$ and $\Real^{n_{\alpha}}$, respectively,
with $n_{\alpha}=|\alpha|$.
Since
\begin{align*}
  B_{\alpha}D^{-1}C_{\alpha}
    & = \frac{1}{k+1} e^{n_{\alpha}} (k+1) \sigma b^* (X_i)^{\sigma-1}  (e^{n_{\alpha}})^T  = \sigma b^* (X_i)^{\sigma-1} e^{n_{\alpha}} (e^{n_{\alpha}})^T,
\end{align*}
it follows that
$A_{\alpha}-B_{\alpha}D^{-1}C_{\alpha}  = \sigma b^* (X_i)^{\sigma-1}
I_{n_{\alpha}} + \sigma(\sigma-1) b^*
(X_i)^{\sigma-2}\usd{\bf x}_{\alpha}(e^{n_{\alpha}})^T,$
which is a sum of a diagonal positive definite matrix and a
$\mathbf{P_0}$-matrix, and thus is a $\mathbf{P}$-matrix.
Therefore, $\det(H_{\alpha})=\det(D)\det(A_{\alpha}-B_{\alpha}D^{-1}C_{\alpha})>0$ for all
$\alpha \subseteq\{1,\ldots,N \}$ with $\alpha\neq\emptyset$, which
implies that $H$ is a $\mathbf{P}$-matrix.\qed

\end{document}